\documentclass[12 pt]{amsart}
\usepackage{amscd,amssymb,amsmath,amsthm}
\input xy
\xyoption{all}
\newtheorem{Proposition}{Proposition}

\newtheorem{Lemma}{Lemma}
\newtheorem{Theorem}{Theorem}

\newtheorem{Conjecture}{Conjecture}

\newcommand{\proj}{\mathbb{P}}
\newcommand{\Z}{\mathbb{Z}}

\newcommand{\rarr}{\rightarrow}

\newcommand{\com}{\mathbb{C}}

\newcommand{\ww}{\widetilde}
\newcommand{\cL}{{\mathcal{L}}}
\newcommand{\q}{\widehat{q}}

\newcommand{\cal}{\mathcal}

\begin{document}

\title{Noether-Lefschetz theory and the
Yau-Zaslow conjecture}
\author{A. Klemm, D. Maulik,  R. Pandharipande and E. Scheidegger}
\date{December 2008}
\maketitle

\begin{abstract}
The Yau-Zaslow conjecture predicts the 
genus 0 curve counts
of $K3$ surfaces in terms of the Dedekind
$\eta$ function. The classical intersection theory of curves 
in the moduli of $K3$ surfaces with Noether-Lefschetz divisors
is related to 3-fold Gromov-Witten invariants via the $K3$ curve counts.
Results by Borcherds and Kudla-Millson determine 
these classical intersections in terms of vector-valued modular
forms. Proven mirror transformations can often be used to calculate the 
3-fold invariants which arise. 

Via a detailed study of the STU model (determining
special curves in the moduli of $K3$ surfaces), we prove the 
Yau-Zaslow conjecture for all curve classes on $K3$ surfaces.
Two modular form identities are required. The first, the Klemm-Lerche-Mayr identity
relating hypergeometric
series to modular forms after mirror transformation, is proven here.
The second, the Harvey-Moore
identity, is proven by D. Zagier and presented in the paper.
\end{abstract}

\tableofcontents

\setcounter{section}{-1}

\baselineskip=18pt

\section{Introduction}

\subsection{Yau-Zaslow conjecture} \label{yzc}
Let $S$ be a nonsingular projective $K3$ surface, and let 
 $$\beta \in \text{Pic}(S) = H^2(S,\mathbb{Z}) \cap H^{1,1}(S,\com)$$
be a nonzero effective curve class.
The moduli space $\overline{M}_{0}(S,\beta)$
of genus 0 stable maps (with no marked points)
has expected dimension
$$\text{dim}^{vir}_\com\left( \overline{M}_{0}(S,\beta) \right)
= \int_\beta c_1(S) + \text{dim}_\com(S) -3 = -1.$$
Hence, the virtual class $[\overline{M}_{0}(S,\beta)]^{vir}$ vanishes,
and the standard Gromov-Witten theory is trivial.

Curve counting on $K3$ surfaces
 is captured instead by the {\em reduced} Gromov-Witten
theory  constructed first via the twistor family in \cite{brl}.
An algebraic construction following \cite{Beh,BehFan} is given in 
\cite{gwnl}. Since the reduced class 
$$[\overline{M}_{0}(S,\beta)]^{red} \in H_0(\overline{M}_{0}(S,\beta), \mathbb{Q})$$
has dimension 0, the reduced Gromov-Witten integrals of $S$,
\begin{equation}\label{veq}
R_{0,\beta}(S) = \int_{[\overline{M}_{0}(S,\beta)]^{red}} 1   
\  \in \mathbb{Q},
\end{equation}
 are
well-defined.
For deformations of $S$ for which $\beta$ remains a $(1,1)$-class,
the integrals \eqref{veq} are invariant.

 The second cohomology of $S$ is a rank 22 lattice
with intersection form 
\begin{equation}\label{ccet}
H^2(S,\mathbb{Z}) \stackrel{\sim}{=} U\oplus U \oplus U \oplus E_8(-1) \oplus E_8(-1)
\end{equation}
where
$$U
= \left( \begin{array}{cc}
0 & 1 \\
1 & 0 \end{array} \right)$$
and 
$$E_8(-1)=  \left( \begin{array}{cccccccc}
 -2&    0 &  1 &   0 &   0 &   0 &   0 & 0\\
    0 &   -2 &   0 &  1 &   0 &   0 &   0 & 0\\
     1 &   0 &   -2 &  1 &   0 &   0 & 0 &  0\\
      0  & 1 &  1 &   -2 &  1 &   0 & 0 & 0\\
      0 &   0 &   0 &  1 &   -2 &  1 & 0&  0\\
      0 &   0&    0 &   0 &  1 &  -2 &  1 & 0\\ 
      0 &   0&    0 &   0 &   0 &  1 &  -2 & 1\\
      0 & 0  & 0 &  0 & 0 & 0 & 1& -2\end{array}\right)$$
is the (negative) Cartan matrix. The intersection form \eqref{ccet}
is even.

The {\em divisibility} $m(\beta)$ is
the maximal positive integer dividing the lattice
element $\beta\in H^2(S,\mathbb{Z})$.
If the divisibility is 1,
$\beta$ is {\em primitive}.
Elements with
equal divisibility and norm  are equivalent up to orthogonal transformation 
 of $H^2(S,\mathbb{Z})$.
By straightforward deformation arguments using the Torelli theorem for
$K3$ surfaces, $R_{0,\beta}(S)$ depends, for effective classes, 
 {\em only} on the divisibility $m(\beta)$ and the norm
$\langle \beta,\beta \rangle$. We will omit the argument $S$ in the
notation.

The genus 0 BPS counts associated to $K3$ surfaces have the following
definition. Let
$\alpha\in \text{Pic}(S)$ be a nonzero class which is both effective and
primitive.
The Gromov-Witten potential $F_{{\alpha}}(v)$ 
for classes proportional
to ${\alpha}$
is 
$${F}_{{\alpha}}=
   \sum_{m>0}\   R_{0,m\alpha}\  v^{m{\alpha}}.$$
The
BPS counts $r_{0,m\alpha}$ are uniquely defined 
by via the Aspinwall-Morrison formula,
\begin{equation}
\label{ccr}
F_\alpha =      \sum_{m>0} \
 r_{0,m\alpha} \ \sum_{d>0}
\frac{v^{dm\alpha}}{d^3}, 
\end{equation}
for both primitive and
divisible classes.

The Yau-Zaslow conjecture \cite{yauz} predicts 
the values of the genus 0 BPS counts for the
reduced Gromov-Witten theory of $K3$ surfaces. We interpret the conjecture in
two parts.

\begin{Conjecture} \label{xxx1}
The BPS count $r_{0,\beta}$ depends upon $\beta$ only through the norm $\langle \beta,\beta \rangle$.
\end{Conjecture}

Conjecture \ref{xxx1} is rather surprising from the point
of view of Gromov-Witten theory since $R_{0,\beta}$
certainly depends upon the divisibility of $\beta$.
Let $r_{0,m,h}$ denote the genus 0 BPS count associated
to a class $\beta$  of divisibility $m$ satisfying
$$\langle \beta,\beta \rangle = 2h-2.$$
Assuming Conjecture \ref{xxx1} holds, we define
$$r_{0,h}= r_{0,m,h}$$
independent{\footnote{Independence of $m$ holds
when $2m^2$ divides $2h-2$. Otherwise, no such class
$\beta$ exists and $r_{0,m,h}$ is defined
to vanish.}} of $m$.

\begin{Conjecture} \label{xxx2}
The BPS counts $r_{0,h}$ are uniquely determined by
\begin{equation} \label{gvre}
 \sum_{h\geq 0} r_{0,h}\ q^h = {\prod_{n=1}^\infty (1 - q^n)^{-24}}.
\end{equation}
\end{Conjecture}
\noindent Conjecture \ref{xxx2} can be written in terms
of the Dedekind $\eta$ function
$$ \sum_{h\geq 0} r_{0,h}\ q^{h-1} =  \eta(\tau)^{-24}$$
 where $q=e^{2\pi i \tau}$.

The conjectures have been previously proven in very few cases.
A proof of the Yau-Zaslow formula \eqref{gvre}
for primitive classes $\beta$ 
via 
Euler characteristics of compactified Jacobians 
following \cite{yauz} 
can be found in \cite{beu,xc,fgd}.
The Yau-Zaslow formula \eqref{gvre}
 was proven via Gromov-Witten theory for
primitive classes $\beta$ by Bryan and Leung \cite{brl}. 
An early calculation by Gathmann \cite{Gat} for a class
$\beta$ of
divisibility 2 was important
for the correct formulation of the conjectures.  
Conjectures 1 and 2 have been proven in the divisibility 2 case
by Lee and Leung \cite{ll} and Wu \cite{wwuu}.
The main result of the paper 
is a proof of Conjectures \ref{xxx1} and \ref{xxx2}
in all cases.

\begin{Theorem} \label{yzz}
The Yau-Zaslow conjecture holds for all nonzero effective classes 
$\beta\in \text{\em Pic}(S)$ on a $K3$ surface $S$.
\end{Theorem}

\subsection{Noether-Lefschetz theory} \label{nnll}
\subsubsection{Lattice polarization}
Let $S$ be a $K3$ surface.
A primitive class 
$L\in \text{Pic}(S)$ is a {\em quasi-polarization}
if
$$\langle L,L \rangle >0  \ \ \ \text{and} \ \  \ \langle L,[C]\rangle
 \geq 0 $$
for every curve $C\subset S$.
A sufficiently high tensor power $L^n$
of a quasi-polarization is base point free and determines
a birational morphism
$$S\rightarrow \tilde{S}$$
contracting A-D-E configurations of $(-2)$-curves on $S$.
Hence, every quasi-polarized $K3$ surface is algebraic.

Let $\Lambda$ be a fixed rank $r$  
primitive{\footnote{An embedding
of lattices is primitive if the quotient is torsion free.}}
embedding
\begin{equation*} 
\Lambda \subset U\oplus U \oplus U \oplus E_8(-1) \oplus E_8(-1)
\end{equation*}
with signature $(1,r-1)$, and
let 
$v_1,\ldots, v_r \in \Lambda$ be an integral basis.
The discriminant is
$$\Delta(\Lambda) = (-1)^{r-1} \det
\begin{pmatrix}
\langle v_{1},v_{1}\rangle & \cdots & \langle v_{1},v_{r}\rangle  \\
\vdots & \ddots & \vdots \\
\langle v_{r},v_{1}\rangle & \cdots & \langle v_{r},v_{r}\rangle 
\end{pmatrix}\ .$$
The sign is chosen so $\Delta(\Lambda)>0$.

A {\em $\Lambda$-polarization} of a $K3$ surface $S$   
is a primitive embedding
$$j: \Lambda \rightarrow \mathrm{Pic}(S)$$ 
satisfying two properties:
\begin{enumerate}
\item[(i)] the lattice pairs 
$\Lambda \subset U^3\oplus E_8(-1)^2$ and
$\Lambda\subset 
H^2(S,\mathbb{Z})$ are isomorphic
 via an isometry which restricts to the identity on $\Lambda$,
 \item[(ii)]
$\text{Im}(j)$ contains
a {quasi-polarization}. 
\end{enumerate}
By (ii), every $\Lambda$-polarized $K3$ surface is algebraic.

The period domain $M$ of Hodge structures of type $(1,20,1)$ on the lattice 
$U^3 \oplus E_8(-1)^2$   is
an analytic open set
 of the 20-dimensional  nonsingular isotropic 
quadric $Q$,
$$M\subset Q\subset \proj\big(    (U^3 \oplus E_8(-1)^2 )    
\otimes_\Z \com\big).$$
Let $M_\Lambda\subset M$ be the locus of vectors orthogonal to 
the entire sublattice $\Lambda \subset U^3 \oplus E_8(-1)^2$.

Let $\Gamma$ be the isometry group of the lattice 
$U^3 \oplus E_8(-1)^2$, and let
 $$\Gamma_\Lambda \subset \Gamma$$ be the
subgroup  restricting to the identity on $\Lambda$.
By global Torelli,
the moduli space $\mathcal{M}_{\Lambda}$ 
of $\Lambda$-polarized $K3$ surfaces 
is the quotient
$$\mathcal{M}_\Lambda = M_\Lambda/\Gamma_\Lambda.$$
We refer the reader to \cite{dolga} for a detailed
discussion.

\subsubsection{Families}

Let $X$ be a 
compact 3-dimensional complex manifold equipped with
holomorphic line bundles
$$L_1, \ldots, L_r  \ \rightarrow X$$
 and a holomorphic map
 $$\pi: X \rarr C$$
to a nonsingular complete curve.

The tuple $(X,L_1,\ldots, L_r, \pi)$ is a
{\em 1-parameter family of nonsingular $\Lambda$-polarized
$K3$ surfaces}
if 
\begin{enumerate}
\item[(i)] the fibers $(X_\xi, L_{1,\xi}, \ldots, L_{r,\xi})$
are $\Lambda$-polarized $K3$ surfaces via
$$v_i \mapsto L_{i,\xi}$$
for every $\xi\in C$,
\item[(ii)] there exists a $\lambda^\pi\in \Lambda$
which is a quasi-polarization of all fibers of $\pi$ 
simultaneously.
\end{enumerate}
The family $\pi$ yields a morphism,
$$\iota_\pi: C \rarr \mathcal{M}_{\Lambda},$$
to the moduli space of $\Lambda$-polarized $K3$ surfaces. 

%Since for each fiber $X_\xi$, the lattice $\Lambda \subset
%\text{Pic}(X_\xi)$ contains a quasi-polarization, we can find a 
%single $\lambda^\pi \in \Lambda$ which is a quasi-polarization
%for all fibers of $\pi$ simultaneously.
Let $\lambda^{\pi}= \lambda^\pi_1 v_1+\dots +\lambda^\pi_r v_r$.
A vector $(d_1,\ldots,d_r)$ of integers is {\em positive} if
$$\sum_{i=1}^r \lambda^\pi_i d_i >0.$$
If $\beta \in \text{Pic}(X_\xi)$ has intersection numbers
$$d_i = \langle L_{i,\xi},\beta \rangle,$$
then $\beta$ has positive degree with respect to the
quasi-polarization if and only if  $(d_1,\dots,d_r)$
is positive.

\subsubsection{Noether-Lefschetz divisors}
Noether-Lefschetz numbers are defined in \cite{gwnl}
by the intersection of $\iota_\pi(C)$ with Noether-Lefschetz 
divisors in $\mathcal{M}_\Lambda$.
We briefly review the definition of the 
Noether-Lefschetz divisors.

Let $(\mathbb{L}, \iota)$ be a rank $r+1$
lattice  $\mathbb{L}$
with an even symmetric bilinear form $\langle,\rangle$ and a primitive embedding
$$\iota: \Lambda \rightarrow \mathbb{L}.$$
Two data sets 
$(\mathbb{L},\iota)$ and $(\mathbb{L}',  \iota')$
are isomorphic if there is an isometry which restricts to identity on $\Lambda$.
The first invariant of the data $(\mathbb{L}, \iota)$ is
the discriminant
 $\Delta \in \mathbb{Z}$ of 
$\mathbb{L}$.

An additional invariant of $(\mathbb{L}, \iota)$ can be 
obtained by considering 
any vector $v\in \mathbb{L}$ for which
\begin{equation}\label{ccff} 
\mathbb{L} = \iota(\Lambda) \oplus \mathbb{Z}v.
\end{equation}
The pairing
$$\langle v, \cdot \rangle : \Lambda \rightarrow \mathbb{Z}$$
determines an element of $\delta_v\in \Lambda^*$.
Let 
$G = \Lambda^{*}/\Lambda$
be quotient defined via the injection
$\Lambda \rightarrow \Lambda^*$
 obtained from the pairing $\langle,\rangle$ on $\Lambda$.
The group $G$ is abelian of order equal to the 
discriminant $\Delta(\Lambda)$.
The image 
$$\delta \in G/\pm$$
of $\delta_v$ is easily seen to be independent of $v$ satisfying 
\eqref{ccff}. The invariant $\delta$ is the {\em coset} of $(\mathbb{L},\iota)$

By elementary arguments, two data sets $(\mathbb{L},\iota)$ and $(\mathbb{L}',\iota')$
of  rank $r+1$ are isomorphic if and only if the discriminants and cosets are
equal.

Let $v_1,\ldots, v_r$ be an integral basis of $\Lambda$ as before.
The pairing of $\mathbb{L}$ 
with respect to an extended basis $v_{1}, \dots, v_{r},v$
is encoded in the matrix
$$\mathbb{L}_{h,d_{1},\dots,d_{r}} = 
\begin{pmatrix}
\langle v_{1},v_{1}\rangle & \cdots & \langle v_{1},v_{r}\rangle & d_{1} \\
\vdots & \ddots & \vdots & \vdots\\
\langle v_{r},v_{1}\rangle & \cdots & \langle v_{r},v_{r}\rangle & d_{r}\\
d_{1} & \cdots & d_{r} & 2h-2
\end{pmatrix}.$$
The discriminant is
$$\Delta(h,d_{1},\dots,d_{r}) 
= (-1)^r\mathrm{det}(\mathbb{L}_{h,d_{1},\dots,d_{r}}).$$
The coset $\delta(h, d_{1},\dots,d_{r})$ is represented by the functional
$$v_i \mapsto d_i.$$

The Noether-Lefschetz divisor $P_{\Delta,\delta} \subset \mathcal{M}_{\Lambda}$
is the closure of the locus of $\Lambda$-polarized $K3$ surfaces $S$ for which
$(\mathrm{Pic}(S),j)$ has rank $r+1$, discriminant $\Delta$, and coset $\delta$.
By the Hodge index theorem, $P_{\Delta,\delta}$ is empty unless $\Delta > 0.$

Let $h, d_{1}, \dots, d_{r}$ determine a positive discriminant
$$\Delta(h,d_{1},\dots,d_{r}) > 0.$$  The Noether-Lefschetz divisor
$D_{h, (d_{1},\dots,d_{r})}\subset \mathcal{M}_{\Lambda}$ is defined by 
the weighted sum
$$D_{h,(d_{1},\dots,d_{r})} 
= \sum_{\Delta,\delta} m(h,d_1,\dots,d_r|\Delta,\delta)\cdot[P_{\Delta,\delta}]$$
where the multiplicity $m(h,d_1,\dots,d_r|\Delta,\delta)$ is the number
of elements $\beta$ of the lattice $(\mathbb{L},\iota)$ of 
type $(\Delta,\delta)$ satisfying
\begin{equation}\label{34f}
\langle \beta, \beta \rangle = 2h-2,\ \  \langle \beta, v_{i}\rangle = d_{i}.
\end{equation}
If the multiplicity is nonzero, then $\Delta | \Delta(h,d_{1},\dots,d_{r})$ so only 
finitely many divisors appear in the 
above sum.

If $\Delta(h,d_{1},\dots,d_{r}) = 0$, the divisor $D_{h,(d_1,\dots,d_r)}$
has an alternate definition.
The tautological
line bundle $\mathcal{O}(-1)$ is $\Gamma$-equivariant on the period domain
$M_\Lambda$ and descends to the {\em Hodge line bundle} 
$$\mathcal{K} \rightarrow \mathcal{M}_{\Lambda}.$$
We define
$D_{h,(d_{1},\dots,d_{r})} = \mathcal{K}^{*}$.
See \cite{gwnl} for an alternate view of degenerate intersection.

If $\Delta(h,d_{1},\dots,d_{r}) < 0$, the divisor 
$D_{h,(d_1,\dots,d_r)}$ on $\mathcal{M}_{\Lambda}$ is defined to vanish
by the Hodge index theorem.

\subsubsection{Noether-Lefschetz numbers}
Let $\Lambda$ be a lattice of discriminant $l=\Delta(\Lambda)$, and
let $(X,L_1,\ldots,L_r,\pi)$ be 
a 1-parameter family of $\Lambda$-polarized $K3$ surfaces.
The Noether-Lefschetz number $NL^\pi_{h,d_1,\dots,d_r}$ is
the  classical intersection
product
\begin{equation}\label{def11}
NL^\pi_{h,(d_1,\dots,d_r)} =\int_C \iota_\pi^*[D_{h,(d_1,\dots,d_r)}].
\end{equation}

Let $\mathrm{Mp}_{2}(\mathbb{Z})$ be
the metaplectic double cover of $SL_{2}(\mathbb{Z})$. 
There is a canonical representation \cite{borch}
associated to $\Lambda$,
$$\rho_{\Lambda}^{*}: 
\mathrm{Mp}_{2}(\mathbb{Z}) \rightarrow \mathrm{End}(\mathbb{C}[G]).$$
The full set of Noether-Lefschetz numbers
$NL^\pi_{h,d_1,\dots,d_r}$ defines a vector valued
modular form
 $$\Phi^{\pi}(q) = \sum_{\gamma\in G} \Phi^{\pi}_{\gamma}(q)v_{\gamma} \in \com[[q^{\frac{1}{2l}}]]
\otimes \com[G],$$
of weight $\frac{22-r}{2}$ and type $\rho_\Lambda^*$
by results{\footnote{While the results of the papers \cite{borch, kudmil}
have considerable overlap, we will follow the point of view of Borcherds.}}
 of Borcherds and Kudla-Millson \cite{borch,kudmil}.
The Noether-Lefschetz numbers are the coefficients{\footnote{If $f$ is a series in $q$,
 $f[k]$ denotes the coefficient of $q^k$.}}
 of the components of 
$\Phi^\pi$,
$$NL^{\pi}_{h,(d_1,\dots,d_r)} = \Phi^{\pi}_{\gamma}\left[ \frac{\Delta(h,d_1,\dots,d_r)}{2l}\right]$$
where $\delta(h,d_1,\dots,d_r) = \pm\gamma$.
The modular form results significantly constrain the Noether-Lefschetz numbers.

\subsubsection{Refinements}
If  $d_1,\ldots,d_r$ do not simultaneously
vanish, refined Noether-Lefschetz divisors
are defined.
If $\Delta(h,d_1,\dots,d_r)>0$, 
$$D_{m,h,(d_1,\dots,d_r)}
\subset D_{h,(d_1,\dots,d_r)}$$ is defined
by requiring the class $\beta \in \text{Pic}(S)$ to satisfy \eqref{34f} and
have divisibility $m>0$. If $\Delta(h,d_1,\dots,d_r)=0$, then
$$D_{m,h,(d_1,\dots,d_r)}=D_{h,(d_1,\dots,d_r)}$$
if $m>0$ is the greatest common divisor of $d_1,\ldots,d_r$
and 0 otherwise.

Refined
Noether-Lefschetz numbers are defined by
\begin{equation}\label{def112}
NL^\pi_{m,h,(d_1,\dots,d_r)} =\int_C \iota_\pi^*[D_{m,h,(d_1,\dots,d_r)}].
\end{equation}
In Section \ref{nl2}, the full 
set of Noether-Lefschetz numbers $NL^\pi_{h,(d_1,\dots,d_r)}$ is
easily shown to determine the refined numbers $NL^\pi_{m,h,(d_1,\dots,d_r)}$.

\subsection{Three theories}
The main geometric idea in the proof  
is the relationship  of three theories associated to   
a 1-parameter family $$\pi:X \rightarrow C$$
 of $\Lambda$-polarized $K3$ surfaces:
\begin{enumerate}
\item[(i)] the Noether-Lefschetz numbers  of $\pi$,
\item[(ii)] the genus 0 Gromov-Witten invariants of $X$,
\item[(iii)] the genus 0 reduced Gromov-Witten invariants of the
$K3$ fibers. 
\end{enumerate}
The Noether-Lefschetz numbers (i) are classical intersection
products while the Gromov-Witten invariants (ii)-(iii) 
are quantum in origin. 
For (ii),
we view the theory in terms  the
Gopakumar-Vafa invariants{\footnote{A review of the definitions
can be found in Section \ref{nl2}.} \cite{GV1,GV2}.

Let $n_{0,(d_1,\dots,d_r)}^X$ denote the Gopakumar-Vafa invariant of $X$
in genus $0$ for $\pi$-vertical curve classes of degrees $d_1,\dots,d_r$
with respect to the line bundles $L_1,\dots, L_r$. Let
$r_{0,m,h}$ denote the reduced $K3$ invariant defined in Section \ref{yzc}.
The following result is proven{\footnote{The result of
the \cite{gwnl} is stated in the rank $r=1$ case, but the
argument is identical for arbitrary $r$.}} in \cite{gwnl}
by a comparison of
the reduced and usual deformation theories of maps of curves
to
the $K3$ fibers of $\pi$.

\begin{Theorem} \label{ffc}
 For degrees $(d_1,\dots,d_r )$ positive with respect to the
quasi-polariza\-tion $\lambda^\pi$,
$$n_{0,(d_1,\dots,d_r)}^X= \sum_{h=0}^\infty \sum_{m=1}^{\infty}
r_{0,m,h}\cdot  NL_{m,h,(d_1,\dots,d_r)}^\pi.$$
\end{Theorem}

\subsection{Proof of Theorem \ref{yzz}}

The STU model described in Section \ref{nl1} is a special family of
rank 2 lattice polarized $K3$ surfaces 
$$\pi^{STU}:X^{STU}\rightarrow \proj^1.$$ 
The fibered $K3$ surfaces of the STU model are themselves
elliptically fibered.
The proof of Theorem \ref{yzz} proceeds in four basic steps:
\begin{enumerate}
\item[(i)] 
The modular form \cite{borch,kudmil} determining the intersections
of the base $\proj^1$ 
with the Noether-Lefschetz divisors is calculated.
For the STU model, the modular form has vector dimension 1 and is
proportional to the product $E_4 E_6$ of Eisenstein series.
\item[(ii)]
Theorem \ref{ffc} is used to show the 3-fold BPS counts
$n_{0,(d_1,d_2)}^{X^{STU}}$ then {\em determine} all the reduced $K3$ 
invariants $r_{0,m,h}$. Strong use is made of the rank 2 lattice of
the STU model. 
\item[(iii)]
The BPS counts $n_{0,(d_1,d_2)}^{X^{STU}}$ are calculated
via mirror symmetry. Since the STU model is realized as a Calabi-Yau
complete intersection in a nonsingular toric variety, the genus 0
Gromov-Witten invariants are obtained after proven mirror transformations
from hypergeometric series. The Klemm-Lerche-Mayr
identity, proven in Section \ref{cfe}, shows the
invariants $n_{0,(d_1,d_2)}^{X^{STU}}$ are themselves related to modular forms.
\item[(iv)] 
Theorem \ref{yzz}  then follows from the Harvey-Moore
identity which simultaneously relates the modular structures of
$$n_{0,(d_1,d_2)}^{X^{STU}}, \ \ r_{0,m,h}, \ \ \text{and}\ \  
NL^{\pi^{STU}}_{m,h,(d_1,d_2)}$$
 in the
form specified by Theorem \ref{ffc}.
D. Zagier's proof of the Harvey-Moore identity is
presented in Section \ref{cew}.
\end{enumerate}

The strategy of proof is special to genus 0. Much less is known in higher
genus.  The Katz-Klemm-Vafa conjecture \cite{kkv,gwnl} for the
integral{\footnote{The integrand $\lambda_g$ is the top Chern class
of the Hodge bundle on $\overline{M}_g(X,\beta)$.}}
$$\int_{[\overline{M}_g(S,\beta)]^{red}} (-1)^g \lambda_g$$
is a particular generalization of 
the Yau-Zaslow formula to higher genera. The KKV formula
does not yet appear easily 
approachable
in Gromov-Witten theory.\footnote{For $g=1$,  
the KKV formula follows for all classes on $K3$ surfaces 
from the Yau-Zaslow formula
via 
the boundary relation for $\lambda_1$.}
However, a proof of the KKV formula for primitive $K3$ classes
in the conjecturally equivalent theory
of stable pairs in the derived category is given
in \cite{ky,pt3}.

 \subsection{Acknowledgments}
We thank R. Borcherds, J. Bruinier, J. Bry\-an, B. Conrad,
 I. Dolgachev, J. Lee,
 E. Looijenga, G. Moore, Y. Ruan, R. Thomas,
G. Tian, W. Zhang, and A. Zinger
 for conversations
about Noether-Lefschetz divisors, reduced invariants of
$K3$ surfaces, and modular forms.
We are especially grateful to D. Zagier for providing us a proof
of the Harvey-Moore identity.

The Clay institute
workshop on $K3$ surfaces in March 2008 which all of us
attended 
played a crucial role in the completion of the paper.
We thank J. Carlson, D. Ellwood, and the Clay institute
staff for providing a wonderful research environment. 
Parts of the paper were written during visits to the
Scuola Normale Superiore in Pisa and the Hausdorff Institute
for Mathematics in Bonn in the summer of 2008.

A.K. was partially supported by DOE grant
DE-FG02-95ER40896.
D.M. was partially supported by a Clay research fellowship.
R.P. was partially support by NSF grant DMS-0500187.
%E.S. was partially supported by ...

\section{The STU model}
\label{nl1}

\subsection{Overview}
The STU model{\footnote{The model has been studied in physics
since the
80's. The
letter $S$ stands for the dilaton and $T$ and $U$
label the torus moduli in the heterotic string. The STU model
was an important example for the
duality between type IIA and heterotic strings
formulated in \cite{kv}. The ideas developed in
 \cite{hm1,hm2,germans,klm,MM} about the STU model
play an important role in our paper.}}
is a particular nonsingular projective
Calabi-Yau 3-fold $X$
equipped with a fibration 
\begin{equation}\label{fibbw}
\pi: X \rightarrow \proj^1.
\end{equation}
Except for 528 points $\xi \in \proj^1$,
the fibers 
$$X_\xi= \pi^{-1}(\xi)$$ are nonsingular
elliptically fibered $K3$ surfaces.
The 528 singular fibers $X_\xi$ have exactly 1 ordinary
double point singularity each.

The 3-fold $X$ is constructed as an anticanonical section of 
a nonsingular projective toric 4-fold $Y$.
The Picard rank of $Y$ is 6. The fibration \eqref{fibbw}
is obtained from a nonsingular toric fibration
$$\pi^Y: Y \rightarrow \proj^1.$$
The image of
$$\text{Pic}(Y)\rightarrow \text{Pic}(X_\xi)$$
determines a rank 2 sublattice of each fiber $\text{Pic}(X_\xi)$
with intersection form
$$
 \left( \begin{array}{cc}
0 & 1 \\
1 & 0 \end{array} \right)\ .$$

The toric data describing the construction of $X\subset Y$ 
and the fibration structure are explained here.

\subsection{Toric varieties}
Let $N$ be a lattice of rank $d$,
$$N\stackrel{\sim}{=} \mathbb{Z}^d.$$
A fan $\Sigma$ in $N$ is a collection of strongly
convex
rational
polyhedral cones 
containing all faces and intersections. A 
toric variety $V_\Sigma$
is
canonically associated to $\Sigma$. The variety $V_\Sigma$ is 
complete of dimension $d$ if the support of $\Sigma$ covers 
$N\otimes_{\mathbb{Z}} \mathbb{R}$.
If all cones are simplicial and if all maximal
cones are generated by a lattice basis, then $V_\Sigma$
is
 nonsingular. See \cite{Cox,Fulton,Oda} for the basic
properties of toric varieties.

Let $\Sigma$ be a fan corresponding to a nonsingular complete
toric variety.
A 1-dimensional cone of $\Sigma$
is a ray with a unique primitive vector.
Let $\Sigma^{(1)}$ denote the set of 1-dimensional cones of $\Sigma$
indexed
by their primitive vectors  
\begin{equation}\label{pgen}
\{ \rho_1, \ldots, \rho_n\}.
\end{equation} 
Let $r^{1}, \ldots r^{\ell}$ be a basis over the integers of
the module of relations among the vectors \eqref{pgen}. We
write the $j^{th}$
relation as 
$$r^j_1 \rho_1 + \ldots + r^j_n \rho_n = 0.$$
Define a torus
  $$(\com^*)^\ell \stackrel{\sim}{=} \prod_{j=1}^\ell \com_j^*$$ 
with factors indexed by the relations.

A simple description
of $V_\Sigma$ is obtained via a quotient construction.
Let $\{z_i\}_{1\leq i \leq n}$ be coordinates on $\com^n$
corresponding to the primitives $\rho_i$ of the rays in
$\Sigma^{(1)}$. 
An action of $\com_j^*$ on $\com^n$ is defined by
\begin{equation}
  \label{eq:weights} 
  \lambda_j\cdot \big(z_1, \dots, z_n \big)   =
  \big(\lambda_j^{r_1^{j}}z_1,\dots,\lambda_j^{r_n^{j}}z_n \big), \ \ \
\lambda_j \in \com_j^*
\end{equation} 
%\footnote{We
%  will use the same symbol $\rho_\ast$ to denote the generators in
%  $\Sigma^{(1)}$ and the corresponding primitive lattice vectors in
%  $N$.} $\rho_i$.  
In order to obtain a well-behaved quotient 
for the induced $(\com^*)^\ell$-action on $\com^n$, 
 an exceptional set $Z(\Sigma)\subset\com^n$
consisting of a finite union of linear subspaces is excluded.
The linear space defined by $\{z_i=0\,| i\in I\}$ is contained
in $Z(\Sigma)$ if there is no single cone in  $\Sigma$
containing all of the primitives $\{ \rho_i\}_{ i \in I}$.
After removing $Z(\Sigma)$, the quotient
\begin{equation}
  \label{eq:toricvariety} 
  V_\Sigma=
    \Big(\com^n \setminus Z(\Sigma)\Big)
    \Big/
     \big(\com^*\big)^\ell 
\end{equation} 
yields the toric variety associated to $\Sigma$.

Since $\ell=n-d$,
the complex dimension of the quotient
 $V_\Sigma$ equals the rank $d$ of the lattice
$N$.
The variety $V_\Sigma$ is equipped with the action of
the quotient torus 
$$T= (\com^*)^n /(\com^*)^\ell.$$
The rank of $\text{Pic}(V_\Sigma)$ is $\ell$.
The primitives $\rho_i$ are in 1--to--1 correspondence with the
$T$-invariant divisors $D_i$ on $V_\Sigma$ 
defined by
\begin{equation}
  D_i = \big\{z_i = 0 \big\} ~\subset V_\Sigma .
\end{equation}
Conversely, the homogeneous coordinate $z_i$ is a section of the line
bundle $\mathcal{O}(D_i)$.
The anticanonical divisor class of $V_\Sigma$ is determined by
\begin{equation} \label{vvq}
-K_{V_\Sigma} = \sum_{i=1}^n D_i.
\end{equation}

\subsection{The toric 4-fold $Y$} \label{YYY}
The fan $\Sigma$ in $\mathbb{Z}^4$ defining the
toric 4-fold $Y$  has 10 rays with primitive elements

\vspace{+10pt}
\begin{center}
\begin{tabular}{lll}
$\rho_1= (1,0,2,3)\ \ \ $ & $\rho_2=(-1,0,2,3)\ \ \ $  & \\
$\rho_3=(0,1,2,3)$ &  $\rho_4=(0,-1,2,3)$  & \\
$\rho_5= (0,0,2,3)$ &  $\rho_6=(0,0,-1,0)$ &  $\rho_7=(0,0,0,-1)$\\
$\rho_8 = (0,0,1,2)$  & $\rho_9=(0,0,0,1)$ &  $\rho_{10}=(0,0,1,1)$.
\end{tabular}
\end{center}
\vspace{+10pt}

\noindent The full fan $\Sigma$ is
obtained from the convex hull of the 10
primitives. 
By explicitly checking each of 24 dimension 4 cones, $Y$ is seen to be a
complete nonsingular toric 4-fold.

Generators $r^1, \ldots, r^6$ of the rank 6 module of relations
among the primitives can be taken to be

\vspace{+10pt}
\begin{center}
\begin{tabular}{rrrrrrrrrrr}
$\rho_1$ & $+\rho_2$ & & & &  $+4\rho_6$ & $+6 \rho_7$ & & & &$=0$\\
&    & $\rho_3$ & $+\rho_4$ & & $+4\rho_6$ & $+6 \rho_7$ & & & &$=0$\\
& & & & $\rho_5$&  $+2\rho_6$ & $+3\rho_7$   & & & &$=0$\\
& & & & & $\rho_6$& $+2 \rho_7$& $+\rho_8$ & & &$=0$\\
& & & & &    & $+\hspace{+5pt} \rho_7$& & $+\rho_9$ &  &$=0$\\
& & & & & $\rho_6$& $+\hspace{+5pt} \rho_7$&  & & $+\rho_{10}$&$=0$\\
\end{tabular}
\end{center}
\vspace{+10pt}

By the identification \eqref{vvq} of $-K_Y$, the product
$\prod_{i=1}^{10} z_i$
defines an anticanonical section. Hence, every product
$$\prod_{i=1}^{10} z_i^{m_i} , \ \ \ m_i\geq 0$$
which is homogeneous of degree $\sum_{i=1}^{10} r^j_i$
with respect to the action \eqref{eq:weights}  of $\com_j^*$
also defines an anticanonical section.
Hence,
\begin{align} 
\label{gt}
{z_{{1}}}^{12}{z_{{4}}}^{12}{z_{{5}}}^{6}{z_{{8}}}^{4}{z_{{9}}}^{2}{z_{{10}}}^{3}, & &
{z_{{1}}}^{12}{z_{{3}}}^{12}{z_{{5}}}^{6}{z_{{8}}}^{4}{z_{{9}}}^{2}{z_{{10}}}^{3},\\
\nonumber
{z_{{2}}}^{12}{z_{{4}}}^{12}{z_{{5}}}^{6}{z_{{8}}}^{4}{z_{{9}}}^{2}{z_{{10}}}^{3}, & &
{z_{{2}}}^{12}{z_{{3}}}^{12}{z_{{5}}}^{6}{z_{{8}}}^{4}{z_{{9}}}^{2}{z_{{10}}}^{3},
\\ \nonumber
{z_{{6}}}^{3}z_{{8}}{z_{{9}}}^{2}, & &
{z_{{7}}}^{2}z_{{10}}
\end{align}
are all sections of $-K_Y$.

From the definitions, we find $Z(\Sigma)$ consists of the union
of the following 11  linear spaces of dimension 2 in $\com^4$,

\vspace{+10pt}
\begin{equation}
  \begin{aligned}
    \label{eq:primitive} 
    I_1 &= \{1,2\}, & I_2 &= \{3,4\}, & I_3 &= \{5,6\}, & I_4 &= \{5,7\}, \\
    I_5 &= \{5,9\}, & I_6 &= \{6,8\}, & I_7 &= \{6,10\}, & I_8 &= \{7,8\},\\
    I_9 &= \{7,9\}, & I_{10} &= \{8,10\}, & I_{11} &= \{9,10\}\ .   
  \end{aligned} 
\end{equation}
\vspace{+10pt}

\noindent Recall, $I_k$ indexes the coordinates which vanish.

A simple verification show
the 6 sections \eqref{gt} of $-K_Y$ do not have a common zero 
on the prequotient $\com^n \setminus Z(\Sigma)$. Hence,
$-K_Y$ is generated by global sections on $Y$. A hypersurface
$$X\subset Y$$
defined by a generic section of $-K_Y$
is nonsingular
by Bertini's Theorem. By adjunction,
$X$ is Calabi-Yau.

\subsection{Fibrations}
The toric variety $Y$ admits two obvious fibrations
$$\pi^Y: Y \rightarrow \proj^1, \ \ \mu^Y: \rightarrow \proj^1$$
given in homogeneous coordinates by
$$\pi^Y( z_1, \ldots, z_{10}) = [z_1,z_2], \ \ \ 
 \mu^Y(z_1, \ldots, z_{10}) = [z_3,z_4].$$
Since $Z(\Sigma)$ contains the linear spaces
$$I_1 = \{1,2\}, \ \ \ I_2 = \{3,4\},$$
both  $\pi^Y$ and $\mu^Y$ are well-defined.

Consider first $\pi^Y$. The fibers of $\pi^Y$ are nonsingular
complete toric 3-folds defined by the fan in 
$$\mathbb{Z}^3 \subset \mathbb{Z}^4, \ \ \ (c_1,c_2,c_3) \mapsto (0,
c_1,c_2,c_3)$$
determined by the primitives $\rho_3 , \ldots, \rho_{10}$.

Let $X$ be obtained from a generic section of $-K_Y$. Let
$$\pi: X \rightarrow \proj^1$$
be the restriction $\pi^Y|_X$.

\begin{Proposition}\label{pww}
Except for 528 points $\xi \in \proj^1$,
the fibers 
$$X_\xi= \pi^{-1}(\xi)$$ are nonsingular
elliptically fibered $K3$ surfaces
The 528 singular fibers $X_\xi$ each have exactly 1 ordinary
double point singularity.
\end{Proposition}

\begin{proof}
Let $P_{k,k}(z_1,z_2|z_3,z_4)$ denote a bihomogeneous polynomial of 
degree $k$ in $(z_1,z_2)$ and degree $k$ in $(z_3,z_4)$.
Let
 $$F=P_{12,12}(z_1,z_2|z_3,z_4),\ \ G=P_{8,8}(z_1,z_2|z_3,z_4), \ \ 
H=P_{4,4}(z_1,z_2|z_3,z_4)$$
be bihomogeneous polynomials.
Then
\begin{equation}\label{pqww}
F z_5^6 z_8^4z_9^2z_{10}^3, \ \ 
G z_5^4z_6z_8^3z_9^2z_{10}^2, \ \ H z_5^2 z_6^2 z_8^2 z_9^2z_{10},
\ \ z_6^3z_8z_9^2, \ \ z_7^2z_{10}
\end{equation}
all determine sections of $-K_Y$.

Let $X$ be defined by a generic linear combination of the
sections \eqref{pqww}. Since the base point free system \eqref{gt}
is contained in \eqref{pqww}, $X$ is nonsingular.
We will prove all the fibers $X_\xi$ are nonsingular, except
for finitely many with exactly 1 ordinary double point each,
 by an explicit study of the equations.

Since $I_7=\{6,10\}$, $I_{10}=\{8,10\}$, and $I_{11}=\{9,10\}$
are in $Z(\Sigma)$, we easily see $X \cap D_{10} = \emptyset$
if the coefficient of $z_6^3z_8z_9^2$ is nonzero.
Similarly
$$X \cap D_{8} = \emptyset, \ \ X \cap D_{9} = \emptyset.$$
Hence, using the last 3 factors of the torus $(\com^*)^\ell$,
the coordinates $z_8$, $z_9$, and $z_{10}$ can all be set to 1.
The equation for $X$ simplifies to
$$Fz_5^6+G z_5^4z_6+ H z_5^2z_6^2+ \alpha z_6^3 + \beta z_7^2.$$

The coordinates $z_1$ and $z_2$ do not simultaneously vanish on $Y$.
There are two charts to consider. By symmetry, the analysis
on each is
identical, so we assume $z_1\neq 0$.
Using the first factor of $(\com^*)^\ell$, we set $z_1=1$.
By the same reasoning, we set $z_3=1$ using the second factor of $(\com^*)^\ell$.
Since $I_3=\{5,6\}$ and $I_4=\{5,7\}$ are in $Z(\Sigma)$ either
$z_5\neq 0 $ or {\em both} $z_6$ and $z_7$ do not vanish.

\vspace{+10pt}
\noindent{\bf Case $z_5\neq 0$.}
Using the third factor of $(\com^*)^\ell$ to set $z_5=1$, we obtain the equation
\begin{equation}\label{ggq}
F(1,z_2|1,z_4)+H(1,z_2|1,z_4) z_6 + G(1,z_2|1,z_4) z_6^2
+ \alpha z_6^3 + \beta z_7^2
\end{equation}
in $\com^4$ with coordinates $z_2,z_4,z_6,z_7$.
The map $\pi$ is given by the $z_2$ coordinate.
The partial derivative of \eqref{ggq} with respect $z_7$ is $2\beta z_7$.
Hence, if $\beta \neq 0$, all singularities of $\pi$ occur when $z_7=0$. 

We need only analyze
the reduced dimension case
\begin{equation}\label{ggq2}
 F(1,z_2|1,z_4)+H(1,z_2|1,z_4) z_6 + G(1,z_2|1,z_4) z_6^2
+ z_6^3 
\end{equation}
with coordinates $z_2,z_4,z_6$.
Here, $\alpha$ has been set to 1 by scaling the equation.
 We must show all the 
fibers of $\pi$  are nonsingular
curves except for finitely many with simple nodes. 
We view equation \eqref{ggq2} as defining a
1-parameter family of paths $\gamma_{z_2}(z_4)$ in the space  
$$\mathcal{C}= \{ \gamma_0 + \gamma_1 z_6 + \gamma_2 z_6^2 + z_6^3\ | \
\gamma_0,\gamma_1,\gamma_2\in \com\}$$ 
of cubic polynomials
in the variable $z_6$. The coordinate of the path is $z_4$.
The variable $z_2$ indexes the family of paths.

Let $\Delta\subset \mathcal{C}$ be the codimension 1 discriminant
locus of cubics with double roots. The discriminant
is irreducible with cuspidal singularities in codimension 2 
in $\mathcal{C}$.
The possible singularities of the fiber $\pi^{-1}(\lambda)$
occur only when the path $\gamma_\lambda(z_4)$
intersects $\Delta$.
The fiber $\pi^{-1}(\lambda)$ is nonsingular over such
an intersection point if either

\begin{enumerate}
\item[(i)] $\gamma_\lambda$ is transverse to $\Delta$
at a nonsingular point of $\Delta$,
\item[(ii)] $\gamma_\lambda$ is transverse to 
the codimension 1 tangent cone of a singular point of $\Delta$.
\end{enumerate}
The fiber $\pi^{-1}(\lambda)$
has a simple node over an intersection point of the
path $\gamma_\lambda(z_4)$ with $\Delta$ if 
\begin{enumerate}
\item[(iii)]
 $\gamma_\lambda$
is tangent to $\Delta$
at a nonsingular point of $\Delta$.
\end{enumerate}
The above are
all the possibilities which can occur in a generic 1-parameter
family of paths in the space of cubic equations.{\footnote{A cusp
of $\pi^{-1}(\lambda)$ occurs, for example, when the path has
contact order 3 at a nonsingular point of the discriminant.}} 
  Possibility
(iii) can happen only for finitely many $\lambda$
and just once for each such $\lambda$.

\vspace{+10pt}
\noindent{\bf Case $z_6\neq 0$ and $z_7\neq 0$}.
Using the third factor of $(\com^*)^\ell$ to set $z_6=1$, we obtain 
the equation
\begin{equation}\label{fvr}
 F(1,z_2|1,z_4)+H(1,z_2|1,z_4) + G(1,z_2|1,z_4)
+ \alpha  + \beta z_7^2
\end{equation}
in $\com^4$ with coordinates
$z_2,z_4,z_5,z_7$.
The partial derivative of \eqref{fvr} with respect $z_7$ is not 0
for $z_7 \neq 0$. Hence, there are no singular fibers of $\pi$ on the chart.
\vspace{+10pt}

We have proven all the fibers $X_\xi$ of $\pi$ are nonsingular except
for finitely many with exactly 1 ordinary double point each.
Let $X_\xi$ be a nonsingular fiber.
Let $$\mu:X \rightarrow \proj^1$$
be the restriction $\mu^Y|_X$.
The fibers of 
product
$$(\pi,\mu): X \rightarrow \proj^1 \times \proj^1$$
are easily seen to be anticanonical sections
of the nonsingular toric surface{\footnote{Since the product 
$(\pi^Y,\mu^Y): Y \rightarrow \proj^1 \times \proj^1$
has fibers isomorphic to the nonsingular complete (hence projective)
toric surface $W$, the 4-fold $Y$ is projective.}}
$W$  with fan in $\mathbb{Z}^2$
determined by the primitives $\rho_5, \ldots, \rho_{10}$.
These anticanonical sections are elliptic curves.
Since $X_\xi$ has trivial canonical bundle by adjunction and
 the map
$$\mu: X_\xi \rightarrow \proj^1$$
is dominant with elliptic fibers, we conclude
$X_\xi$ is an elliptically fibered $K3$ surface.

The Euler characteristic of $X$ can be calculated by 
toric intersection in $Y$,
$$\chi_{top}(X)= -480.$$
The Euler characteristic of a nonsingular $K3$ fibration
over $\proj^1$ is $48$. Since each fiber singularity reduces the
Euler characteristic by 1,
we conclude $\pi$ has exactly 528 singular fibers.
\end{proof}

For emphasis, we will sometimes denote the STU model
by
$$\pi^{STU}: X^{STU} \rightarrow \proj^1.$$

\subsection{Divisor restrictions} \label{YYY2}
The divisors $D_1$, $D_2$, $D_8$, $D_9$, and $D_{10}$ 
have already been shown to restrict to the trivial class in
$\text{Pic}(X_\xi)$.
The divisors $D_3$ and $D_4$ restrict to the fiber class
$F\in \text{Pic}(X_\xi)$ of the
elliptic fibration
\begin{equation}\label{xxdd}
\mu: X_\xi \rightarrow \proj^1.
\end{equation}
Certainly $F^2=0$.
Let $S\in \text{Pic}(X_\xi)$ denote
the restriction of $D_5$. Toric calculations yield
the products
$$F\cdot S =1, \ \ S\cdot S =-2.$$
Hence, $S$ may be viewed as the section class of the
elliptic fibration \eqref{xxdd}.
The divisors $D_6$ and $D_7$ restrict to classes in
the rank 2 lattice generated by $F$ and $S$.

The restriction of $\text{Pic}(Y)$ to each 
fiber $X_\xi$ is a rank 2 lattice generated by
$F$ and $S$ with intersection form
$$
 \left( \begin{array}{cc}
0 & 1 \\
1 & -2 \end{array} \right)\ .$$
We may also choose generators $L_1=F$ and $L_2=F+S$
with intersection form
$$ \Lambda=
 \left( \begin{array}{cc}
0 & 1 \\
1 & 0 \end{array} \right)\ .$$

\subsection{1-parameter families}
\label{nl11}
Let $X$ be a 
compact 3-dimensional complex manifold  equipped with two holomorphic
line bundles
$$L_1,\  L_2 \rightarrow X$$
and a holomorphic
map
$$\pi:X \rarr C$$
to a nonsingular complete curve.

The data $(X,L_1,L_2,\pi)$ determine
 a {\em family of $\Lambda$-polarized $K3$ surfaces} if
the fibers $(X_\xi, L_{1,\xi}, L_{2,\xi})$
are $K3$ surfaces with intersection form 
$$
 \left( \begin{array}{cc}
L_{1,\xi}\cdot L_{1,\xi} & L_{2,\xi}\cdot L_{1,\xi} \\
L_{1,\xi}\cdot L_{2,\xi} & L_{2,\xi} \cdot L_{2,\xi} \end{array} \right)
=\left( \begin{array}{cc}
0 & 1 \\
1 & 0 \end{array} \right)$$
and there exists a simultaneous quasi-polarization.
The 1-parameter family $(X,L_1,L_2,\pi)$ yields a morphism,
$$\iota_\pi: C \rarr \mathcal{M}_{\Lambda},$$
to the moduli space of $\Lambda$-polarized $K3$ surfaces.

The construction $(X^{STU}, L_1,L_2, \pi^{STU})$ 
of the STU model in  Sections \ref{YYY}-\ref{YYY2} is almost
a $1$-parameter
family of $\Lambda$-polarized $K3$ surfaces. The only failing
is the 528 singular fibers of $\pi^{STU}$.
Let
$$\epsilon:C \stackrel{2-1}{\longrightarrow} \proj^1$$
be a hyperelliptic curve branched over the 528 points of
$\proj^1$ corresponding to the singular fibers of $\pi$.
The family
$$\epsilon^*(X^{STU}) \rightarrow C$$
has 3-fold double point singularities over the 528 nodes
of the fibers of the original family.
Let
$$\widetilde{\pi}^{STU}: \widetilde{X}^{STU} \rightarrow C$$
be obtained from a small resolution
$$\widetilde{X}^{STU} \rightarrow \epsilon^*(X^{STU}).$$
Let $\widetilde{L}_i\rightarrow \widetilde{X}^{STU}$ be the
pull-back of $L_i$ by $\epsilon$.
The data $$(\widetilde{X}^{STU},  \widetilde{L}_1,
\widetilde{L}_2,\widetilde{\pi}^{STU})$$ determine a 1-parameter
family of $\Lambda$-polarized
$K3$ surfaces, see Section 5.3 of \cite{gwnl}.
The simultaneous quasi-polarization is obtained from
the projectivity of $X^{STU}$.

%For every fiber $\ww{X}^{STU}_\xi$, the class $\ww{L}_{1,\xi}$ is nef. Hence,
%$$\langle \ww{L}_{1,\xi}, [C] \rangle \geq 0$$
%for every effective curve $C\subset \ww{X}^{STU}_\xi$.

\subsection{Gromov-Witten invariants}\label{gggr}
Since $X^{STU}$ is defined by an anticanonical section in a semi-positive
nonsingular toric variety $Y$, the genus 0 Gromov-Witten
invariants have been proven by Givental \cite{giv1, giv2,lly,pgiv} 
to be related 
by mirror transformation
to
hypergeometric solutions of the Picard-Fuchs equations of
the Batyrev-Borisov mirror.
By Section 5.3 of \cite{gwnl}, the Gromov-Witten invariants
of $\widetilde{X}^{STU}$ are exactly twice the
Gromov-Witten invariants of $X^{STU}$ for curve classes in the fibers.

\section{Noether-Lefschetz numbers and reduced $K3$ invariants} \label{nlred}
\subsection{Refined Noether-Lefschetz numbers}\label{reff}
Following the notation of Section \ref{nnll}, let
\begin{equation*} 
\Lambda \subset U\oplus U \oplus U \oplus E_8(-1) \oplus E_8(-1)
\end{equation*}
be primitively embedded with signature $(1,r-1)$ and
integral basis $v_1,\dots,v_r$. Let
$(X,L_1,\ldots,L_r, \pi)$ be a 1-parameter family of
$\Lambda$-polarized $K3$ surfaces.
Let $d_1,\dots,d_r$ be integers which do not all vanish.

\begin{Lemma}
\label{refinednl}
The Noether-Lefschetz numbers
$NL^{\pi}_{h, (d_{1},\dots,d_{r})}$
completely determine the refinements 
$NL^{\pi}_{m, h, (d_{1},\dots,d_{r})}.$
\end{Lemma}
\begin{proof}
By definition, 
the refined Noether-Lefschetz numbers satisfy two elementary identities.
The first is
$$NL^{\pi}_{h,(d_{1},\dots,d_{r})} = 
\sum_{m=1}^{\infty} NL^{\pi}_{m, h, (d_{1},\dots,d_{r})}.$$
If $m$ does not divide all $d_i$, then $NL^{\pi}_{m,h,(d_1,\dots,d_r)}$
vanishes. If $m$ divides all $d_i$, then a second
identity holds:
$$NL^{\pi}_{m,h, (d_{1},\dots,d_{r})} = 
NL^{\pi}_{1, h', (d_{1}/m, \dots, d_{r}/m)}$$ 
where $2h-2 = m^{2}(2h'-2)$.

If $\Delta(h,d_1,\dots,d_r)=0$, the refined number
$NL^\pi_{m,h,(d_1,\dots,d_r)}$ vanishes by definition unless
$m$ is the GCD of $(d_1,\dots,d_r)$. In the latter
case, 
$$NL^\pi_{h,(d_1,\dots,d_r)}=NL^\pi_{m,h,(d_1,\dots,d_r)}.$$
Hence the Lemma is trivial in the $\Delta(h,d_1,\dots,d_r)=0$
case.

If $\Delta(h,d_1,\dots,d_r)>0$,
we prove the
Lemma  by induction on
$\Delta$.  
The second identity reduces us to the case where $m=1$.
The first identity determines the $m=1$ case in terms of 
the Noether-Lefschetz number
$NL_{h,(d_1,\dots,d_r)}$  and refined
numbers  with 
$$\Delta(h',d'_1,\dots,d'_r) < \Delta(h,d_1,\dots,d_r).$$
\end{proof}

\subsection{STU model}
The resolved version of the STU model
$$\widetilde{\pi}^{STU}:\widetilde{X}^{STU}\rightarrow C$$
is  lattice polarized with respect to
$$\Lambda =  \left( \begin{array}{cc}
0 & 1  \\
1 & 0  \end{array} \right).$$
The application of the  results  of \cite{borch,kudmil}
to the STU model is extremely simple.  
Since the lattice $\Lambda$ is unimodular, 
the corresponding representation $\rho_{\Lambda}^{*}$ is 1-dimensional and, 
in fact, is the trivial representation of 
$\mathrm{Mp}_{2}(\mathbb{Z})$.  The Noether-Lefschetz degrees are thus encoded
by a scalar modular form of weight $\frac{22-r}{2} = 10$.  
The space of such forms is well-known to be of  dimension 1 and
spanned by the product of Eisenstein series{\footnote{
The Eisenstein series $E_{2k}$ is the modular form
defined by the equation 
$$-\frac{B_{2k}}{4k} E_{2k}(q) = -\frac{B_{2k}}{4k} + \sum_{n\geq 1} 
\sigma_{2k-1}(n) q^{n},$$
where $B_{2n}$ is the $2n^{th}$ Bernoulli number and
$\sigma_n(k)$ is the sum of the $k^{th}$ powers of
the divisors of $n$,
$$\sigma_k(n) = \sum_{i|n} i^k.$$}}
$$E_{10}(q) = E_{4}(q)E_{6}(q) = 1 - 264 \sum_{n\geq 1} \sigma_{9}(n) q^{n}.$$
Hence, a single Noether-Lefschetz 
calculation determines the full series.  

\begin{Lemma} $NL^{\widetilde{\pi}}_{0,(0,0)} = 1056$.
\end{Lemma}
\begin{proof}
By Proposition \ref{pww}, the STU model 
$$\pi^{STU}: X^{STU} \rightarrow \proj^1$$
has 528 nodal fibers.
Let $S$ be a 
 fiber of the resolved family
 $\widetilde{\pi}^{STU}$  lying over a singular fiber of $\pi$.
The Picard lattice of $S$ certainly contains
\begin{equation} \label{kvt}
 \left( \begin{array}{rrr}
0 & 1 & 0  \\
1 & 0 &0 \\
0 & 0 & -2 \end{array} \right)
\end{equation}
spanned by $L_1$, $L_2$, and the $(-2)$-curve $E$ of the small resolution.
Let 
$$\widetilde{\iota}: C \rightarrow \mathcal{M}_\Lambda$$
be the map to moduli.
Since a class $\beta$ satisfying
$$\langle \beta,\beta\rangle =-2$$
on a $K3$ surface is either effective or anti-effective,
the set theoretic intersections of
$\widetilde{\iota}$ with $D_{0,(0,0)}$ correspond to fibers
of $\widetilde{\pi}$ where $L_1$ and $L_2$ do not generate an
ample class --- precisely the 528 fibers of $\widetilde{\pi}$
lying over the
singular fibers of $\pi$.

The divisor $D_{0,(0,0)}$ has multiplicity exactly 2 at the 528
intersections with $\widetilde{\iota}$
since $E$ and $-E$ are the only $-2$ classes orthogonal to
$L_1$ and $L_2$. Finally, since $E$ has normal bundle $(-1,-1)$
in $\widetilde{X}^{STU}$, the curve $\widetilde{\iota}$ is transverse
to the reduced divisor $\frac{1}{2} D_{0,(0,0)}$ at the 528 intersections.
We conclude $NL^{\widetilde{\pi}}_{0,(0,0)}= 528\cdot 2= 1056$.
\end{proof}

\begin{Proposition}\label{studegrees}
The Noether-Lefschetz degrees of the 
resolved STU model are given by the equation
$$NL^{\tilde{\pi}}_{h, (d_{1},d_{2})} 
= -4E_{4}(q)E_{6}(q) \left[\frac{\Delta(h,d_{1},d_{2})}{2}\right].$$
\end{Proposition}

\subsection{BPS states}
Let 
$(\widetilde{X}^{STU},\ww{L}_1,\ww{L}_2,\widetilde{\pi}^{STU})$ 
the $\Lambda$-polarized STU model
The vertical classes are the kernel of the push-forward map
by $\widetilde{\pi}$,
$$0 \rightarrow H_2(\widetilde{X},\mathbb{Z})^{\widetilde \pi} \rightarrow 
H_2(\widetilde{X},\mathbb{Z})
\rightarrow H_2(C,\mathbb{Z}) \rightarrow 0.$$
While $\widetilde{X}$ need not be a projective variety,
$\widetilde{X}$ carries a $(1,1)$-form $\omega_K$ which is K\"ahler on the
$K3$ fibers of $\widetilde{\pi}$. The existence of a fiberwise K\"ahler
form is sufficient to define Gromov-Witten theory for
vertical classes  
$$0\neq \gamma \in H_2(\widetilde{X},\mathbb{Z})^{\widetilde \pi}.$$
The fiberwise K\"ahler form $\omega_K$ is obtained by a small
perturbation of the quasi-K\"ahler form obtained from the
quasi-polarization. The associated Gromov-Witten theory is
independent of the perturbation used.
%{\footnote{There are
%several alternate ways to define the Gromov-Witten theory for 
%the resolved STU model
%$\widetilde{X}$. 
%The easiest is via symplectic geometry, since by \cite{liruan},
%small resolutions are symplectic. Since $\widetilde{X}$
%is so close to projective, an algebraic approach is also easily found.}}

Let $\overline{M}_{0}(\widetilde{X},\gamma)$ be the moduli space of
stable maps from connected genus $0$ curves to $\widetilde{X}$.
Gromov-Witten theory is defined
by
integration against the virtual class,
\begin{equation}
\label{klk}
N_{0,\gamma}^{\widetilde{X}}
 = \int_{[\overline{M}_{0}(\widetilde{X},\gamma)]^{vir}} 1.
\end{equation}
The expected dimension of the moduli space is 0.

The genus 0 
Gromov-Witten potential $F^{\widetilde{X}}(v)$ for nonzero vertical classes
is the series
$${F}^{\widetilde{X}}=
  \sum_{0\neq \gamma\in H_2({\widetilde{X}},\mathbb{Z})^{\widetilde \pi}}  
 N^{\widetilde{X}}_{0,\gamma} \  v^\gamma$$
where $v$ is the curve class variable.
The
BPS counts $n_{0,\gamma}^{\ww{X}}$
of Gopakumar and Vafa are uniquely defined 
by the following equation:
\begin{equation*}
F^{\widetilde{X}}  =     \sum_{0\neq \gamma\in
H_2(\widetilde{X},\mathbb{Z})^{\widetilde \pi}} 
 n_{0,\gamma}^{\widetilde{X}} \ \sum_{d>0}
\frac{v^{d\gamma}}{d^3}. 
\end{equation*}
Conjecturally, the invariants $n_{0,\gamma}^{\widetilde{X}}$ are integral and
obtained from the cohomology of an as yet unspecified moduli
space of sheaves on $\widetilde{X}$. 
We do not assume the conjectural properties
hold.

Using the  $\Lambda$-polarization,
we define the BPS counts
\begin{equation}
n_{0,(d_1,d_2)}^{\ww{X}} 
= \sum_{\gamma\in H_2(\ww{X},\mathbb{Z})^{\widetilde \pi},
 \ \int_\gamma \ww{L}_i
=d_i} n_{0,\gamma}^{\ww{X}}
\end{equation}
when $d_1$ and $d_2$ are not both 0. 
%The sum has only finitely many nonzero terms.
%{\footnote{For the
%resolved model $\widetilde{X}$, the finiteness claim is
%a consequence of the transition formula of \cite{liruan}.}}
%Since $\ww{L}_1$ is nef on the fibers, $n_{0,(d_1,d_2)}$
%vanishes unless $d_1 \geq 0$.
 
The original STU model,
$$\pi^{STU}:X^{STU} \rightarrow \proj^1,$$ 
with 528 singular fibers
is a nonsingular, projective, Calabi-Yau 3-fold. Hence the  
Gromov-Witten invariants are
well-defined. Let $n_{0,(d_1,d_2)}^X$ denote the
fiberwise Gopakumar-Vafa invariant with degrees $d_i$ measured by $L_i$.
By the argument of Section \ref{gggr},
$$n_{0,(d_1,d_2)}^{\ww{X}} = 2 n_{0,(d_1,d_2)}^X $$
when $d_1$ and $d_2$ are not both 0.

\subsection{Invertibility of constraints}
Let $\mathcal{P}\subset \mathbb{Z}^2$ be the set of pairs 
$$\mathcal{P} = \{ \ (d_1,d_2)\neq (0,0)\ | \ 
 d_1\geq 0, \  d_1 \geq -d_2\ \}\ .$$
Pairs $(d_2,d_2) \in \mathcal{P}$
are certainly positive with respect to any 
quasi-polarization for $\ww{\pi}^{STU}$
since such $(d_1,d_2)$ can be realized by linear combinations
of the effective classes $F$ and $S$.
 
Theorem \ref{ffc}
applied to the resolved STU model yields the equation
\begin{equation}\label{maintheoremforstu}
n^{\widetilde{X}}_{0, (d_1, d_2)} = 
\sum_{h=0}^\infty \sum_{m=1}^{\infty} 
r_{0,m,h}\cdot NL^{\tilde{\pi}}_{m,h,(d_1,d_2)}
\end{equation}
for $(d_1,d_2)\in \mathcal{P}$.
The BPS states on the left side will be computed by mirror
symmetry in Section \ref{cfe}.
The refined Noether-Lefschetz degrees are determined by 
Lemma \ref{refinednl} and
Proposition \ref{studegrees}.
Consequently, 
equation \eqref{maintheoremforstu} provides 
constraints on the reduced $K3$ invariants $r_{0,m,h}$

The integrals $r_{0,m,h}$ are very simple in case $h\leq 0$.
By Lemma 2 of \cite{gwnl},
$r_{0,m,h}=0$ for $h<0$,
$$r_{0,1,0}= 1,$$ 
and $r_{0,m,0}=0$ otherwise.

\begin{Proposition}\label{bbh}
The set of integrals $\{r_{0,m,h}\}_{m\geq 1, h> 0}$ 
is uniquely determined by the set of  constraints \eqref{maintheoremforstu}
 for $(d_1\geq 0,\ d_2> 0)$ and the integrals $r_{0,m,h\leq 0}$.
\end{Proposition}

\begin{proof}
A certain subset of the 
linear equations with $d_2> 0$ will be shown to be
upper triangular in the variables $r_{0,m,h}$.  
Picard rank $2$ is crucial for the argument.

Let us fix in advance the values of $m\geq 1$ and  $h>0$.
We proceed by induction on $m$  assuming the reduced invariants $r_{0,m',h}$
have already been determined for all $m' < m$.  
The assumption is vacuous when $m=1$.
We can also assume 
$r_{0,m,h'}$ has been determined inductively for $ h' < h$.
If $2h-2$ is not divisible by $2m^{2}$, then
we have $r_{0,m,h}= 0$, so we can further assume
$$2h-2 = m^{2}(2s-2)$$
for an integer $s> 0$.

Consider equation \eqref{maintheoremforstu} for 
$(d_1,d_2) = (m(s-1),m)$.  Certainly 
$$NL^{\ww{\pi}}_{m',h',(m(s-1),m)} = 0$$
unless $m'$ divides $m$.  By the Hodge index theorem,
we must have
\begin{equation}\label{fred}
\Delta(h', m(s-1),m) = 2 - 2h' + m^2(2s-2) \geq 0
\end{equation}
if $NL^{\ww{\pi}}_{m,h',(m(s-1),m)} \neq 0$.
Inequality \eqref{fred} implies
$h' \leq h$.

Therefore, the constraint \eqref{maintheoremforstu}
takes the form
$$n^{\widetilde{X}}_{0,(m(s-1),m)} =  
r_{0,m,h} NL^{\widetilde{\pi}}_{m, h, (m(s-1),m)}+ \dots,$$
where the dots represent terms involving $r_{0,m',h'}$ with either 
$$m' < m  \ \ \ \text{ or } \ \ \  m'=m,\ h' < h.$$
The leading coefficient is given by
$$NL^{\widetilde{\pi}}_{m, h, (m(s-1),m)} = 
NL^{\widetilde{\pi}}_{h, (m(s-1),m)} = -4.$$
As the system is upper-triangular, we can invert to solve for $r_{0,m,h}$.
\end{proof}

\subsection{Proof of the Yau-Zaslow conjecture}
By Proposition \ref{bbh},
we need only show the answer for $r_{0,m,h}$
predicted by the Yau-Zaslow conjecture
satisfies the 
 constraints \eqref{maintheoremforstu} 
for all pairs $(d_1\geq 0,\ d_2>0)$.

Let $X^{STU}$ be the original Calabi-Yau 3-fold of the
STU model. 
%Nonzero
%effective $\pi$-vertical classes $(d_1,d_2)$ must be of the
%form
%\begin{equation}\label{bert}
%(d_1> 0, d_2\in \mathbb{Z}) \cup (d_1=0, d_2 > 0).
%\end{equation}
%Since $L_1$ is nef, $d_1\geq 0$. If $d_1=0$, the effective
%class must be a multiple of the fiber $F$ implying $d_2>0$.
Let 
\begin{equation}\label{llqh}
D_2^3 F^X= \sum_{(d_1,d_2)\in \mathcal{P}} 
d_2^{3}\ N^{X}_{0,(d_1,d_2)}\ q_{1}^{d_1}q_{2}^{d_2} 
\end{equation}
be the third derivative{\footnote{$D_2=q_2 \frac{d}{dq_2}$.}}
 of the genus 0 Gromov-Witten series
for  $\pi$-vertical classes in $\mathcal{P}$.

We can calculate $D_2^3 F^X$ by the constraint \eqref{maintheoremforstu}
{\em assuming the validity of the Yau-Zaslow conjecture},
\begin{equation}\label{p45}
D_2^3 F^X = \sum_{(d_1,d_2)\in \mathcal{P}}
 d_2^3\ c(d_1,d_2)\ 
 \frac{q_{1}^{d_1}q_{2}^{d_2}}{1 - q_{1}^{d_1}q_{2}^{d_2}}
\end{equation}
where
$c(k,l)$
is the coefficient of $q^{kl}$ in 
$$-2\frac{E_{4}(q)E_{6}(q)}{\eta^{24}(q)}.$$

\begin{Proposition} \label{ppx}
The Yau-Zaslow conjecture is implied by the identity 
$$\sum_{(d_1,d_2)\in \mathcal{P}} 
d_2^{3}N^{X}_{0,(d_1,d_2)}q_{1}^{d_1}q_{2}^{d_2} 
=\sum_{(d_1,d_2)\in \mathcal{P}}
 d_2^3\ c(d_1,d_2)\ 
 \frac{q_{1}^{d_1}q_{2}^{d_2}}{1 - q_{1}^{d_1}q_{2}^{d_2}}.$$
\end{Proposition}

\begin{proof}
The $q_1^{d_1}q_2^{d_2}$ coefficient of the above identity is simply $d_2^3$
times the constraint \eqref{maintheoremforstu}. Since we
only require the constraints in 
 case $$(d_1\geq 0,\ d_2>0)\in \mathcal{P},$$ the
identity implies all the constraints we need.
\end{proof}

The remainder of the paper is devoted to the proof of 
Proposition \ref{ppx}. 
The genus 0 Gromov-Witten invariants of $X$  are related, after
mirror transformation, to hypergeometric solutions of the associated 
Picard-Fuchs system of 
differential equations. Hence, Proposition \ref{ppx}
amounts to a subtle identity among special functions. 

\label{nl2}

\section{Mirror transform}
\label{cfe}
\subsection{Picard-Fuchs}
Let $\pi:X \rightarrow \proj^1$ be the STU model.
Let $$\delta_0\in H^*(X,\com)$$ denote the identity class.
A basis of $H^2(X,\com)$ is obtained
from the restriction of the
toric divisors of $Y$ discussed in Section \ref{YYY2},
$$\delta_1=2D_1+2D_3+D_5, \ \ \delta_2=D_3, \ \ \delta_3=D_1.$$
Recall, $\delta_3$ vanishes on the fibers of $\pi$.
Let $\{\delta_j\}$ be a full basis of $H^*(X,\com)$ extending
the above selections.

Let 
$u_1,u_2,u_3$ be the
 canonical
coordinates for the mirror family
with respect to the divisor basis $\delta_1,\delta_2,\delta_3$.
Let
$$\theta_i = u_i \frac{\partial}{\partial u_i}.$$
The Picard-Fuchs system associated to the mirror of $X^{STU}$ is:
\begin{equation}
  \begin{aligned}
    \cL_1 &=\theta_{{1}} \left(
      \theta_{{1}}-2\,\theta_{{2}}-2\,\theta_{{3}} \right) -12\,\left(
      6\,\theta_{{1}}-5 \right) \left( 6\,\theta_{{1}}-1 \right)
    u_{{1}}\\
    \cL_2 &= {\theta_{{2}}}^{2}- \left( 2\,\theta_{{2}}+2\,
      \theta_{{3}}-\theta_{{1}}-2 \right) \left(
      2\,\theta_{{2}}+2\,\theta_
      {{3}}-\theta_{{1}}-1 \right) u_{{2}},\\
    \cL_3 &={\theta_{{3}}}^{2}- \left(
      2\,\theta_{{2}}+2\,\theta_{{3}}-\theta_{{1}} -2\right) \left(
      2\,\theta_{{2}}+2\,\theta_{{3}}-\theta_{{1}}-1
    \right) u_{{3}}\,.
  \end{aligned}
  \label{eq:pfF0}
\end{equation}
The system is obtained canonically from the 
Batyrev-Borisov construction, see \cite{CKatz} for the formalism.

\subsection{Solutions}
%Let $T_i$ be graded homogeneous coordinates on $H^*(X,\mathbb{\com})$
%and with $T_0$ corresponding to the identity class and
%$T_1, T_2, T_3$ corresponding to the divisors
%$D_3, D_5, D_1$ respectively.
%The higher $T_i$ correspond to cohomology classes
%in higher codimension.

A fundamental solution to the Picard-Fuchs system can be written
in terms of GKZ hypergeometric series, 
\begin{equation}
\label{gttn}
\varpi\in H^*(X,\mathbb{\com}) \otimes_{\mathbb{\com}} 
\mathbb{\com}[
\log(u_1),\log(u_2),\log(u_3)][[u_1,u_2,u_3]].
\end{equation}
Let
$\varpi(u,\delta_j)$ be the corresponding coefficient
%{\footnote
%{The coefficient is well-defined with respect to 
%a basis of $H^*(X,\com)$.}} 
of \eqref{gttn}, then
$$\cL_i\  \varpi(u,\delta_j)=0.$$
The standard normalization of $\varpi$ satisfies two
important properties:
\begin{enumerate}
\item[(i)]
The $\delta_0$ coefficient is the unique solution
$$\varpi(u,\delta_0)=1 + O(u)$$
holomorphic at $u=0$.
\item[(ii)] For $1\leq i \leq 3$,
$$\varpi(u,\delta_i)= \frac{\varpi(u,\delta_0)}{2\pi i} \log(u_i) + O(u)$$
are the logarithmic solutions.
\end{enumerate}

Let $T_1,T_2,T_3$ be coordinates on $H^2(X,\com)$ with
respect to the basis $\delta$.
% and 
%let
%$$Q_i = \exp(2\pi i T_i).$$
The mirror transformation is
defined by
$$T_i = \frac{\varpi(u,\delta_i)}{\varpi(u,\delta_0)} = \frac{1}{2\pi i}
\log(u_i) + O(u)$$
for $1\leq i \leq 3$.

The mirror transformation relates the genus 0 Gromov-Witten
theory of $X$ to the Picard-Fuchs system for the mirror family.
For anticanonical hypersurfaces in toric varieties, a proof
is given in \cite{giv2}.

%The cone of effective curve classes of
%$X$ has 3 generators 
%$$L^1, L^2,L^3\in H_2(X,\mathbb{Z}).$$ 
%We have already seen that the divisors $D_8$, $D_9$, and $D_{10}$
%of $Y$ do not intersect $X$. We describe $L^i$ by the
%vector of intersections with the divisors 
%$D_j$,
%$$L^i = (l^i_1, \ldots, l^i_7) , \ \ l^i_j = \int_{L^i} D_j.$$
%Also, let
%$$l^i_0 = \int_{L^i}  D_1+ \ldots +D_{10} =  -\int_{L^i} K_Y.$$
%A standard toric analysis yields
%\begin{equation}
%  \label{eq:MoriXt}
%  \renewcommand{\arraystretch}{1.3}
%  \begin{array}{r@{(}r@{,~}r@{,~}r@{,~}r@{,~}r@{,~}r@{,~}r@{)}l} 
%    L^{1}= &1&1&0&0&-2&0&0   \\
%    L^{2}= &0&0&1&1&-2&0&0   \\
%    L^{3}= &0&0&0&0&1&2&3   & 
%  \end{array}\ \ \  ,
%\end{equation}
%and 
%$$l^1_0 =0, \ \ \l^2_0=0, \ \ l^3_0 = -6.$$

\subsection{Mirror transform for $q_3=0$}
We introduce two modular parameters
\begin{equation}
  \begin{aligned}
    \tau_1&= T_1, & \tau_2&= T_1+T_2 \,.
%&\tau_3&= T_3
\\
  \end{aligned}
  \label{eq:defTU}
\end{equation}
For $i=1$ and 2, let
$$\q_i = \exp(2\pi i \tau_i),$$
 and let
$q_3= \exp(2\pi i T_3)$.

Our first
step is to find a modular expression for the mirror map and 
the period $\varpi(u,\delta_0)$ to leading order in $q_3$. 
We prove two formulas discovered by Klemm, Lerche,
and Mayr in  \cite{klm}.

\begin{Lemma} We have \label{vvh}
\begin{equation*}
  \begin{aligned}
  u_1 &=  \frac{2(j(\q_1)+j(\q_2)-\mu)}{j(\q_1)j(\q_2)
  +\sqrt{j(\q_1)(j(\q_1)-\mu)}\sqrt{j(\q_2)(j(\q_2)-\mu)}}+O(q_3), \\
  u_2 &= \frac{(j(\q_1)j(\q_2) +\sqrt{j(\q_1)(j(\q_2)-\mu)}
\sqrt{j(\q_2)(j(\q_2)-\mu)})^2}{{4j(\q_1)j(\q_2)(j(\q_1)+j(\q_2)-\mu)^2}} +O(q_3),
  \end{aligned}
%  \label{eq:jrelation}
\end{equation*}
where $\mu=1728$ and
%\begin{equation}
%  j(q)=\frac{E_4^3}{\Delta}=\frac{1}{q}+744+196884q+O(q^2)\ 
%\end{equation} 
\begin{equation}
  j(q)=\frac{E_4^3}{\eta^{24}}=\frac{1}{q}+744+196884q+O(q^2)\ 
\end{equation} 
is the normalized $j$ function.
\end{Lemma}

\begin{Lemma}\label{vvhh}
$  \text{\em Lim}_{q_3\rightarrow 0}\ 
 \varpi(u,\delta_0)= E_4(\q_1)^{\frac{1}{4}} E_4(\q_2)^{\frac{1}{4}} \, .$
 \end{Lemma}

\begin{proof} We prove Lemmas \ref{vvh} and \ref{vvhh} together.
The first step is to 
perform the following change of variables
\begin{align*}
%  \label{eq:F2toF0}
  u_1 &= z_1, & u_2 &= \frac{z_2}{2}\left(1+\sqrt{1-4z_3}\right), & u_3 &=
  \frac{z_2}{2}\left(1-\sqrt{1-4z_3}\right), 
\end{align*}
{with the inverse change}
\begin{equation*}
%  \label{eq:F2toF0}
  z_1 = u_1,\ \ \ \  z_2 = u_2 + u_3, \ \ \ \
  z_3 = \frac{u_2u_3}{(u_2+u_3)^2}.
\end{equation*}
In the new variables, the limit $u_3\to 0$ becomes the limit
$z_3\to0$.

The statement of Lemma \ref{vvh} in the variables $z_i$
 remains unchanged to first order in $q_3$. We will prove
\begin{equation*}
  \begin{aligned}
  z_1 &= \frac{2\left(j(\q_1)+j(\q_2)-\mu\right)}{j(\q_1)j(\q_2)
  +\sqrt{j(\q_1)(j(\q_1)-\mu)}\sqrt{j(\q_2)(j(\q_2)-\mu)}}+O(q_3),\\
  z_2 &= \frac{(j(\q_1)j(\q_2) +\sqrt{j(\q_1)(j(\q_2)-\mu)}
\sqrt{j(\q_2)(j(\q_2)-\mu)})^2}{{4j(\q_1)j(\q_2)(j(\q_1)+j(\q_2)-\mu)^2}} +O(q_3)\,.
  \end{aligned}
%  \label{eq:jrelationF2}
\end{equation*}

The Picard-Fuchs differential operators~(\ref{eq:pfF0}) can be rewritten as
\begin{equation*}
  \label{eq:pfF2toF0}
  \begin{aligned}
    \cL_1'(z) &= \cL_1(u), \\
z_2\sqrt{1-4z_3} \   \cL_2'(z) &= \cL_2(u) - \cL_3(u),\\
%    \cL_3'(z) &= u_2\theta_3(u)-u_3\theta_2(u),
z_2\sqrt{1-4z_3}\ \cL_3'(z) &= u_3\cL_2(u)-u_2\cL_3(u),
  \end{aligned}
\end{equation*}
with
\begin{equation*}
  \begin{aligned}
    \cL_1'&=\theta_{{1}} \left( \theta_{{1}}-2\,\theta_{{2}}
 \right) -12\, \left( 6\,\theta_{{1}}-5 \right)  \left( 6\,\theta_{{1}}-1 \right) z_{{1}},\\
    \cL_2' &=\theta_{{2}} \left( \theta_{{2}}-2\,\theta_{{3}} \right) - \left( 2\,\theta_{{2}}-\theta_{{1}}-2\right)  \left( 2\,\theta_{{2}}-\theta_{{1}}-1
 \right) z_{{2}},\\
    \cL_3' &={\theta_{{3}}}^{2}- \left( 2\,\theta_{{3}}-\theta_{{2}}-2 \right) 
 \left( 2\,\theta_{{3}}-\theta_{{2}}-1 \right) z_{{3}}
  \end{aligned}
  \label{eq:pfF2}
\end{equation*}
where now $\theta_i=z_i\frac{d}{dz_i}$. 
Since $\cL_3'(z) \to 0$
in the limit $z_3\to 0$, we need only focus on
$\cL_1'(z)$ and $\cL_2'(z)$. 

Next, we transform $\cL_1'(z)$ and $\cL_2'(z)$
to new variables $y_1, y_2 ,y_3$ via the change
\begin{equation*}
  \begin{aligned}
  z_1 &= \frac{2\left(y_1+y_2-\mu\right)}{y_1y_2
  +\sqrt{y_1 (y_1-\mu)}\sqrt{y_2(y_2-\mu)}},\\
  z_2 &=\frac{(y_1y_2 +\sqrt{y_1(y_1-\mu)}\sqrt{y_2(y_2-\mu)})^2}
  {4y_1y_2(y_1+y_2-\mu)^2},\\ 
  z_3 &=y_3.
  \end{aligned}
\end{equation*}
We obtain
\begin{align*}
  \cL''_1&=
    y_1^2 y_2(y_1- \mu)\partial^2_{y_1}
  + y_1y_2(y_1- \frac{\mu}{2})\partial_{y_1}
  - y_1 y_2^2(y_2- \mu) \partial^2_{y_2}
\\ &\phantom{=} 
  - y_1y_2(y_2-\frac{\mu}{2})\partial_{y_2} 
  + 60(y_1-y_2), \\\notag
  \cL''_2&=
  - y_1^2(y_1- \mu)\,\partial^2_{y_1}
  + y_1( \frac{\mu}{2}- y_1)\partial_{y_1}
  + y_2^2(y_2- \mu)\,\partial^2_{y_2}
  + y_2(y_2-\frac{\mu}{2})\partial_{y_2}
  \\ &\phantom{=} 
  - 2 y_1 y_3 (y_1 -\mu) \partial_{y_1} \partial_{y_3}
  + 2 y_2 y_3 (y_2- \mu) \partial_{y_2} \partial_{y_3}.
%   \cL''_3&=
%     (y_1 - y_2)(y_1-\mu)^2\,y_1^2 \partial^2_{y_1} 
%   + 2( y_1-\mu)y_1( y_1^2 + \mu y_2 - 2\,y_1\,y_2) \,\partial_{y_1}\\
%   &\phantom{=} 
%   - 2( y_1-\mu) \,y_1 ( y_1 - y_2)( y_2 -\mu) \,y_2\,\partial_{y_1}\partial_{y_2}
%   + (y_1 - y_2)(y_2 -\mu)^2 y_2^2\partial^2_{y_2}\\
%   &\phantom{=} 
%   - 2(y_2-\mu)y_2(y_2^2 + \mu y_1  -2 y_1 y_2 ) \partial_{y_2}\\ 
%   &\phantom{=} 
%   + ( y_1 - y_2)^3\,z( 4\,z -1) \partial^2_{z}
%     + ( y_1 - y_2)^3 (6\,z -1) \partial_{z}\\
%   &\phantom{=}   
%   - 4 ( y_1 - y_2)^2\,( y_2-\mu)y_2\,z \partial_{y_2}\partial_{z}
%   + 4 (y_1-\mu )y_1 ( y_1 - y_2)^2\,z \partial_{y_1}\partial_{z}.
\end{align*}
In the limit $y_3\to 0$, the second line on the right for $\cL''_2$ vanishes.
We can combine $\cL_1''$ and $\cL_2''$ to obtain the following 
simple forms:
\begin{equation*}
  \begin{aligned}
%    \left( y_1 - y_2 \right){\cal L}'''_1&:= 
{\cal L}_1''+y_1
    \lim_{y_3\rightarrow 0}{\cal L}''_2= \left( y_1 - y_2 \right)
    \,\left( 60 - \left( y_1-\frac{\mu}{2} \right)
      \,y_1\,\partial_{y_1} -
      \left( y_1-\mu \right) \,y_1^2\,\partial_{y_1}^2 \right),\\
%    \left( y_1 - y_2 \right){\cal L}'''_2&:= 
{\cal L}_1''+y_2
    \lim_{y_3\rightarrow 0}{\cal L}''_2= \left( y_1 - y_2 \right)
    \,\left( 60 - \left( y_2-\frac{\mu}{2} \right)
      \,y_2\,\partial_{y_2} -
      \left( y_2-\mu \right) \,y_2^2\,\partial_{y_2}^2 \right).\\
  \end{aligned}
%  \label{eq:pf3}
\end{equation*}

The solution $\varpi(y,\delta_0)_{y_3=0}$ therefore 
satisfies the differential equation
% factorizes in the leading order in $u_3$ into two
%functions $g_1(y_1)$ and $g_2(y_2)$ both fulfilling the 
%differential equation  
\begin{equation} 
  {\cal L}=\left(y-\mu \right) y^2\partial_y^2+
  \left( y-\frac{\mu}{2} \right) y\,\partial_{y} -  60  \ .
  \label{eq:diffj}
\end{equation}
in both $y_1$ and $y_2$.

Changing \eqref{eq:diffj} to the variable $t=\frac{1728}{y}$ yields
\begin{equation*} 
  {\cal L}=
t(1-t) \partial^2_t+ (1 - \tfrac{3}{2} t)\partial_t-\tfrac{5}{144}\ , 
\end{equation*} 
which by comparing with the general hypergeometric differential 
operator 
\begin{equation*}
  {\cal L}=t(1-t) \partial^2_t+ (c - (1+a+b)t)\partial_t-a b
\end{equation*} 
is identified with the system 
$$_2F_1(a,b;c;t)={}_2F_1(\frac{1}{12},\frac{5}{12};1;t(\tau)).$$ 
According to the results
 of Klein and Fricke as reviewed in ~\cite{zag1}, 
 we have a unique (up to scaling) solution $g_0$ to (\ref{eq:diffj})
locally analytic at $y=\infty$. The solution
 can be written as
\begin{equation*}
  g_0(j(\tau))=(E_4)^{\frac{1}{4}}(\tau), 
\qquad \qquad  y(\tau)=j(\tau) \ .
\end{equation*} 
Moreover, the inverse is
$$\tau(y)=\frac{g_1(y)}{2 \pi i g_0(y)},$$
where $g_1$ is a logarithmic solution at $y=\infty$ of 
${\cal L}$, unique up to normalization and 
addition of $g_0$.

Transformation of the solution $\varpi(u,\delta_0)$ is seen 
to be analytic in a neighborhood of $t_1=t_2=0$. We conclude
$$\varpi(u,\delta_0)_{u_3=0}= E_4^{\frac{1}{4}}(\tau_1)
E_4^{\frac{1}{4}}(\tau_2).$$
By comparing the first few coefficients of the actual solutions 
$\varpi(u,\delta_i)$  in the $u_3\rightarrow 0$ limit,
 we can uniquely identify 
\begin{equation*}
  \tau_1(u)= T_1(u),\qquad \tau_2(u)= T_1(u)+T_2(u)\ .  
\end{equation*}
Hence, Lemma \ref{vvhh} is
 established. Lemma \ref{vvh} is proven by transforming  
back to the $u_1$ and $u_2$ variables.
\end{proof}

Restricted to a $K3$ fiber of $\pi: X \rarr \proj^1$, we
have
$$\delta_1 = 2F+S, \ \ \delta_2 = F.$$
The coordinates $2\pi i\tau_1$ and $2\pi i\tau_2$ correspond to the
divisor basis
$$L_2=F+S, \ \ L_1=F$$
of the $K3$ fiber. Since the variables $q_1$ and $q_2$ of Section \ref{nlred}
measure degrees against $L_1$ and $L_2$, we see
$$\q_1=q_2 \ \  \text{and} \ \  \q_2=q_1$$
for the fiber geometry.

\label{mttt}
\subsection{B-model}

The mirror transformation results of Section \ref{mttt}
together with a B-model calculation of the periods will
be used to prove the following result discovered by
Klemm, Mayr, and Lerche \cite{klm}.

\begin{Proposition} \label{btxg} We have
$$2+
\sum_{(d_1,d_2)\in \mathcal{P}} 
d_2^{3}N^{X}_{0,(d_1,d_2)}q_{1}^{d_1}q_{2}^{d_2} 
=
2
\frac{E_4(q_1)E_6(q_1)}{\eta^{24}(q_1)} 
\frac{E_4(q_2)}{j(q_1)-j(q_2)}.$$
\end{Proposition}

The left side  of Proposition \ref{btxg}
is the left side of Proposition \ref{ppx}
with an added
degree 0 constant 2.

\begin{proof}
We will use following universal
expression for the Gromov-Witten invariants
of $X$ in terms of the periods of the mirror:
\begin{multline*}2+
\sum_{(d_1,d_2)\in \mathcal{P}} 
d_2^{3}N^{X}_{0,(d_1,d_2)}q_{1}^{d_1}q_{2}^{d_2} = \\
\lim_{q_3\rightarrow 0} \ 
\frac{1}{\varpi(u(T),\delta_0)^2}
\sum_{i,j,k=1}^3 \frac{\partial u_i}{\partial \tau_1}
\frac{\partial u_j}{\partial \tau_1} \frac{\partial u_k}{\partial \tau_1}
Y_{i,j,k}(u(T))
\end{multline*}
where $Y_{i,j,k}$ are the Yukawa couplings of the mirror
family, see \cite{CKatz,klm}.
%We have used here that $T_2$ corresponds to $F+S$ on the
%$K3$ fibers of $\pi$.

The periods $Y_{i,j,k}$ can be explicitly computed 
via Griffith transversality \cite{klm} and greatly 
simplify in the $q_3 \rightarrow 0$ limit.
We tabulate the results below:
\begin{equation*}
  \begin{aligned}
    Y_{111} &= \frac{8(1-\,\tilde u_{{1}})}{\tilde u_{{1}}^3\Delta_1},
    & Y_{133} &= \frac{2\tilde u_{{1}}(1-\tilde u_{{1}})}{\tilde
      u_{{3}}\Delta_1},\\ 
    Y_{112} &= \frac{2(1-\tilde u_{{1}})^2+\tilde u_{{1}}^2(\tilde
      u_{{2}}-\tilde u_{{3}})}{\tilde u_{{1}}^2\tilde u_{{2}}\Delta_1},
    & Y_{222} &= \frac{\left(1-2\tilde u_{{1}}\right)A_2}{2\tilde
      u_{{2}}^2\Delta_1\Delta_2},\\
    Y_{113} &= \frac{2(1-\tilde u_{{1}})^2+\tilde u_{{1}}^2(\tilde
      u_{{3}}-\tilde u_{{2}})}{\tilde u_{{1}}^2\tilde
      u_{{3}}\Delta_1},
    & Y_{223} &= \frac{\left(1-2\tilde u_{{1}}\right)A_3}{2\tilde
      u_{{3}}\tilde u_{{2}}\Delta_1\Delta_2},\\
    Y_{122} &= \frac{2\tilde u_{{1}}(1-\tilde u_{{1}})}{\tilde
      u_{{2}}\Delta_1}, 
    & Y_{233} &= \frac{\left(1-2\tilde u_{{1}}\right)A_2}{2\tilde
      u_{{3}}\tilde u_{{2}}\Delta_1\Delta_2},\\
    Y_{123} &= \frac{(1-\tilde u_{{1}})\left((1-\tilde
        u_{{1}})^2-(\tilde u_{{2}}+\tilde u_{{3}})\tilde
        u_{{1}}^2\right)}{\tilde u_{{1}}\tilde u_{{2}}\tilde
      u_{{3}}\Delta_1},
    & Y_{333} &= \frac{\left(1-2\tilde u_{{1}}\right)A_3}{2\tilde
      u_{{3}}^2\Delta_1\Delta_2}.
  \end{aligned}
%  \label{eq:byukawa}
\end{equation*}
Here, we have 
 introduced the variables
$$\tilde u_1=432 u_1, \ \ \tilde u_2=4 u_2, \ \ \tilde u_3=4 u_3$$ and 
the discriminant loci 
\begin{equation}
  \label{eq:5}
  \begin{aligned}
    \Delta_1 &= (1-\tilde u_{{1}})^4 -2(\tilde u_{{2}}+\tilde u_{{3}})\tilde u_{{1}}^2(1-\tilde u_{{1}})^2 + (\tilde u_{{2}}-\tilde u_{{3}})^2\tilde u_{{1}}^4,\\
    \Delta_2 &= (1-\tilde u_{{2}}-\tilde u_{{3}})^2-4\tilde u_{{2}}\tilde u_{{3}}.
  \end{aligned}
\end{equation}
The quantities $A_2$ and $A_3$ are defined by
\begin{equation}
  \label{eq:37}
  \begin{aligned}
    A_2 &= \left( 1+\tilde u_{{2}}-\tilde u_{{3}} \right)  \left( 1-\tilde u_{{1}} \right) ^{2}+{\tilde u_{{
1}}}^{2} \left( 1-\tilde u_{{3}}-3\,\tilde u_{{2}} \right)  \left( \tilde u_{{2}}-\tilde u_{{3}}
 \right), \\
    A_3 &= \left( 1+\tilde u_{{3}}-\tilde u_{{2}} \right)  \left( 1-\tilde u_{{1}} \right) ^{2}+{\tilde u_{{
1}}}^{2} \left( 1-\tilde u_{{2}}-3\,\tilde u_{{3}} \right)  \left( \tilde u_{{3}}-\tilde u_{{2}}
 \right).
  \end{aligned}
\end{equation}
The normalizations of the Yukawa couplings
 $Y_{i,j,k}$  are fixed by the classical 
intersections. 

The leading behavior of the mirror map for $u_1,u_2$ is obtained 
by rewriting Lemma \ref{vvh}
 in terms of $E_4(\tau_i)$ and $E_6(\tau_i)$ as
\begin{equation} 
  \begin{aligned}
  u_1&=\frac{1}{864}\left(1 - \frac{E_6(\tau_1)\,E_6(\tau_2)}{{E_4(\tau_1)}^{\frac{3}{2}}\,
{E_4(\tau_2)}^{\frac{3}{2}}}\right)+{\cal O}(q_3)\ ,  \\
  u_2&=\frac{\left( {E_4(\tau_1)}^3 - {E_6(\tau_1)}^2 \right) \,\left( {E_4(\tau_2)}^3 - 
{E_6(\tau_2)}^2 \right) }{4\,{\left( {E_4(\tau_1)}^{\frac{3}{2}}\,{E_4(\tau_2)}^{\frac{3}{2}} 
- E_6(\tau_1)\,E_6(\tau_2) \right) }^2}+{\cal O}(q_3)
  \end{aligned}
\end{equation}
Denote the leading behavior of the last mirror map by
\begin{equation} 
  u_3=q_3 f_3(\q_1,\q_2)+{\cal O}(q_3^2)\ .
\end{equation}
The derivatives of the mirror maps with respect to $T_2$ are easily 
evaluated using the standard identities
\begin{equation*}
  \begin{aligned}
   q \frac{d}{dq} E_2&=\tfrac{1}{12}(E_2^2-E_4)\\ 
    q\frac{d}{dq} E_4&=\tfrac{1}{3}(E_2 E_4-E_6)\\  
    q \frac{d}{dq} E_6&=\tfrac{1}{2}(E_2 E_6-E_4^2)\\  
    %j&=\frac{E_4^3}{E_4^3-E_6^2}, 
    q\frac{d}{dq} j&= -j \frac{E_6}{E_4}.
  \end{aligned}
\end{equation*} 
We find, to leading order in $q_3$,
\begin{equation*}
\begin{array}{rl}
\frac{\partial u_1}{\partial \tau_1}&=\frac{E_6(\tau_2)\,\left( {E_4(\tau_1)}^3 - {E_6(\tau_1)}^2 \right) }
  {1728\,{E_4(\tau_2)}^{\frac{3}{2}}\,{E_4(\tau_1)}^{\frac{5}{2}}}\\
\frac{\partial u_2}{\partial \tau_1}
&=\frac{{\sqrt{E_4(\tau_1)}}\,\left( {E_4(\tau_2)}^3 - 
{E_6(\tau_2)}^2 \right) \,
    \left( -\left( {E_4(\tau_1)}^{\frac{3}{2}}\,E_6(\tau_2) \right)  + 
      {E_4(\tau_2)}^{\frac{3}{2}}\,E_6(\tau_1) \right) \,\left( {E_4(\tau_1)}^3 - {E_6(\tau_1)}^2 \right) 
    }{4\,{\left( {E_4(\tau_2)}^{\frac{3}{2}}\,{E_4(\tau_1)}^{\frac{3}{2}} - E_6(\tau_2)\,E_6(\tau_1)
        \right) }^3}
\end{array}
\end{equation*}
The derivative 
$\frac{\partial u_3}{\partial \tau_1}$ can be  written to this order as
\begin{equation} 
  \frac{\partial u_3}{\partial \tau_1}
= \frac{u_3}{f_3(\q_1,\q_2)}\frac{\partial}{\partial
\tau_1} f_3(\q_1,\q_2)+O(u_3^2)\ .
  \label{eq:dzdt}
\end{equation}

There are many simplifications in the limit $u_3\rightarrow 0$. First 
the triple couplings 
$$Y_{133}, \ \ Y_{233}, \ \ Y_{333}$$  do not have enough inverse 
powers of $u_3$ and therefore do
 not contribute by the vanishing (\ref{eq:dzdt}). 
Second, the surviving $Y_{i,j,k}$ simplify in the limit. 
We evaluate  
\begin{multline} \label{v678}
 \lim_{q_3\rightarrow 0}
\frac{1}{\varpi(u(T),\delta_0)^2}
\sum_{i,j,k=1}^3 \frac{\partial u_i}{\partial \tau_1}
\frac{\partial u_j}{\partial \tau_1} \frac{\partial u_k}{\partial \tau_1}
Y_{i,j,k}(u(T))
=\\
  -2\,\frac{{E_4(\tau_2)\,E_4(\tau_1)}\,E_6(\tau_2)\,
    \left( {E_4(\tau_1)}^3 - {E_6(\tau_1)}^2 \right) }{ {E_4(\tau_2)}^3\,{E_6(\tau_1)}^2-{E_4(\tau_1)}^3\,{E_6(\tau_2)}^2}\ .
\end{multline}
The possible linear 
dependence on $f_3(\q_1,\q_2)$ drops out as claimed in \cite{klm}!
Using the standard identities
$$j= \frac{E_4^3}{\eta^{24}}, \ \ \eta^{24}= E_4^3-E_6^2,$$
we obtain the right side of Proposition \ref{btxg}.
\end{proof}

\section{The Harvey-Moore identity}
\label{cew}
\subsection{Proof of Proposition \ref{ppx}}
After evaluating the left side via Proposition \ref{btxg}
and dividing by 2,
Proposition \ref{ppx} amounts to a modular form
identity.
%Let 
%$$S_k \subset M_k \subset M^!_k$$ denote the spaces
%of cusp form, modular forms, and weakly holomorphic{\footnote{
%Holomorphic
%except for a possible pole at infinity.}}  modular 
%forms for $\text{SL}(2,\mathbb{Z})$.
Let
$$f(\tau)= \frac{E_4(\tau)E_6(\tau)}{\eta(\tau)^{24}} =
\sum_{n=-1}^\infty c(n) q^n $$
where $q=\exp(2\pi i \tau)$.
Then, we must prove
\begin{equation}
\label{harmoo}
\frac{f(\tau_1) E_4(\tau_2)}{j(\tau_1)-j(\tau_2)}
= \frac{q_1}{q_1-q_2} + E_4(\tau_2) - \sum_{d,k,\ell>0}
\ell^3 c(k\ell)\ q_1^{kd}q_2^{\ell d}.
\end{equation}
Equation \eqref{harmoo} is the Harvey-Moore identity
conjectured in \cite{hm1}.

\subsection{Zagier's proof of the Harvey-Moore identity}
The Harvey-Moore identity implies
Proposition \ref{ppx} and concludes the
proof of the Yau-Zaslow
conjecture. We present here Zagier's argument
from \cite{zag2}.

Let 
$S_k \subset M_k \subset M^!_k$ denote the spaces
of cusp forms, modular forms, and weakly holomorphic
{\footnote{Holomorphic
except for a possible pole at infinity.}}  modular 
forms for $\Gamma=\text{SL}(2,\mathbb{Z})$.
Certainly
$$f(\tau) \in M_{-2}^!.$$

For each $n\geq 0$, there is a unique function $F_n\in M_4^!$
satisfying
$$F_n(\tau)= q^{-n} + \mathcal{O}(q)$$
as $\mathfrak{I}(\tau) \rightarrow \infty$.
Uniqueness follows from the vanishing of $S_4$.
Existence follows by writing $F_n(\tau)$ as $E_4(\tau)$ times
a polynomial in $j(\tau)$,
$$F_0= E_4, \ \ F_1= E_4(j-984), \ \ F_2= 
E_4(j^2-1728 j+393768)\ \ldots \ .$$
We draw several consequences:
\begin{enumerate}
\item[(i)]
$F_1|T_n=n^3 F_n$ for all $n\geq 1$, where $T_n$
is the $n^{th}$ Hecke operator in weight 4.
Indeed, $T_n$ sends $M_4^!$ to itself and, by standard
formulas for the action of $T_n$ on Fourier expansions,
$T_n$ sends $q^{-1}+ \mathcal{O}(q)$ to 
$n^3 q^{-n}+ \mathcal{O}(q)$.
\item[(ii)] $F_1= - f'''$ where prime
denotes differentiation by 
$$\frac{1}{2\pi i}\frac{d}{d\tau} = q \frac{d}{dq}.$$
We see $f'''$ lies in $M_4^!$  by the $k=4$
case of Bol's identity
$$\frac{d^{k-1}}{d\tau^{k-1}} (f|_{2-k} \gamma) =
\left( \frac{d^{k-1} f}{d\tau^{k-1}}\right)|_k \gamma \ \ \ 
\forall \gamma\in \Gamma.$$
Since, the Fourier expansion of $f'''$ begins
as $-q^{-1} + \mathcal{O}(q)$, the claim is proven.
\item[(iii)] For $\mathfrak{I}(\tau_1) > \text{max}_{\gamma\in \Gamma}\  
\mathfrak{I}(\gamma\tau_2)$,
\begin{equation*}
\label{bllp}
\frac{ f(\tau_1) E_4(\tau_2)}{j(\tau_1)-j(\tau_2)} =
\sum_{n=0}^\infty F_n(\tau_2) q_1^n.
\end{equation*}
Let $L(\tau_1,\tau_2)$ denote the
 left side of \eqref{bllp}. We see
$L(\tau_1,\tau_2)$ is a meromorphic modular
form in $\tau_2$ with a simple pole of residue
$-\frac{1}{2\pi i}$ at $\tau_2=\tau_1$ (since $j'=-E_4^2 E_6/\eta^{24}$)
and no poles outside $\Gamma \tau_1$. Moreover, $L(\tau_1,\tau_2)$
tends to 0 as $\mathfrak{I}(\tau_2) \rightarrow \infty$.
These properties characterize $L(\tau_1,\tau_2)$ uniquely
and show that the $n^{th}$ Fourier coefficient with respect to
$\tau_1$ for $\mathfrak{I}(\tau_1)\rightarrow \infty$ has
the properties characterizing $F_n(\tau_2)$.
\end{enumerate}

Combining (i) and (ii)
with the formula for the action of $T_n$ on Fourier expansions,
we obtain,
\begin{eqnarray} \label{vp2}
F_n(\tau) & = & (-n^{-3} f''')|T_n = n^{-3}\left(
q^{-1} - \sum_{m=1}^\infty m^3 c(m)\ q^m \right)| T_n\\ \nonumber
& = & q^{-n} - \sum_{\stackrel{k,\ell, d >0}{kd=n}} \ell^3 c(k\ell)\  
q^{\ell d}\ 
\end{eqnarray}
for $n>0$.
The Harvey-Moore identity follows from
 \eqref{vp2} and (iii) together
with the equality $F_0=E_4$. \qed

\pagebreak

\vspace{+7pt}
\noindent Departments of Physics \\
\noindent  University of Bonn and University of Wisconsin

\vspace{+7pt}
\noindent Department of Mathematics\\
\noindent Columbia University

\vspace{+7pt}
\noindent Department of Mathematics\\
\noindent Princeton University

\vspace{+7pt}
\noindent Department of Mathematics\\
\noindent  University of Augsburg

\end{document}